\documentclass[MSNbibl,nameyear,seceqn,dvips]{arxstspdf}
\usepackage{flushend}
\usepackage{stfloats}
\usepackage{graphicx}

\volume{26}
\issue{1}
\pubyear{2011}
\firstpage{21}
\lastpage{46}
\doi{10.1214/10-STS345}

\makeatletter

\newtheorem{Theorem}{Theorem}[section]
\newproclaim{Example}[Theorem]{Example}

\newcommand{\E}{\mathbf{E}}
\newcommand{\EE}{\mathbf{E}}
\newcommand{\Cov}{\operatorname{Cov}}
\newcommand{\Var}{\operatorname{Var}}
\newcommand{\N}{\mathcal{N}}

\fnbelowfloat

\setattribute{abstract}   {width}  {36pc}
\setattribute{keyword}    {width}  {36pc}

\makeatother

\begin{document}
\begin{frontmatter}

\title{Feature Matching in Time Series Modeling}
\relateddoi{T1}{Discussed in \doi{10.1214/11-STS345A},
\doi{10.1214/11-STS345C}, \doi{10.1214/11-STS345B} and \doi
{10.1214/11-STS345D}; rejoinder at \doi{10.1214/11-STS345REJ}.}
\runtitle{Feature Matching}

\begin{aug}
\author{\fnms{Yingcun} \snm{Xia}\ead[label=e1]{staxyc@nus.edu.sg}}
\and
\author{\fnms{Howell} \snm{Tong}\ead[label=e2]{howell.tong@gmail.com}}
\runauthor{Y. Xia and H. Tong}

\affiliation{National University of Singapore, London School of
Economics and National University of Singapore, London School of
Economics}

\address{Yingcun Xia is Professor of Statistics, Department of Statistics and Applied
Probability, Risk  Management Institute, National University of Singapore, Singapore \printead{e1}.
Howell Tong is Emeritus Chair Professor of Statistics, London School of Economics,
Houghton Street, London WC2A 2AE,
United Kingdom \printead{e2}.}

\end{aug}

%
\begin{abstract}
Using a time series model to mimic an observed time series has a long
history. However, with regard to this objective, conventional
estimation methods for discrete-time dynamical models are frequently
found to be wanting. In fact, they are characteristically misguided in
at least two respects: (i) assuming that there is a true model; (ii)
evaluating the efficacy of the estimation as if the postulated model
is true. There are numerous examples of models, when fitted by
conventional methods, that fail to capture some of the most basic
global features of the data, such as cycles with good matching periods,
singularities of spectral density functions (especially at the origin)
and others. We argue that the shortcomings need not always be due to
the model formulation but the inadequacy of the conventional fitting
methods. After all, all models are wrong, but some are useful \textit{if
they are fitted properly}. The practical issue becomes one of how to
best fit the model to data.

Thus, in the absence of a true model, we prefer an alternative
approach to conventional model fitting that typically involves
one-step-ahead prediction errors. Our primary aim is to match the joint
probability distribution of the observable time series, including
long-term features of the dynamics that underpin the data, such as
cycles, long memory and others, rather than short-term prediction. For
want of a better name, we call this specific aim \textit{feature
matching}.

The challenges of model misspecification, measurement errors and the
scarcity of data are forever present in real time series modeling. In
this paper, by synthesizing earlier attempts into an
extended-likelihood, we develop a~systematic approach to empirical time
series analysis to address these challenges and to aim at achieving
better feature matching. Rigorous proofs are included but relegated to
the \hyperref[appendix]{Appendix}. Numerical results, based on both
simulations and real
data, suggest that the proposed
catch-all approach has several advantages over the conventional
methods, especially when the time series is short or with strong
cyclical fluctuations. We conclude with listing directions that require
further development.
\end{abstract}

%

\end{frontmatter}

\renewcommand{\theequation}{\arabic{section}.\arabic{equation}}
\section{Introduction}

Dynamical models, either in continuous time or in discrete time,
have been widely used to describe the changing world. Interestingly,
salient features of many seemingly complex observations can
sometimes be captured by simple dynamical models, as demonstrated
most eloquently by Sir Isaac Newton in the seventeenth century when
he used his model, Newton's law of universal gravitation, to explain
Kepler's observations concerning planetary motion. In statistics,
dynamical models are the \textit{raison d'{\^{e}}tre} of time series
analysis. For a time series, the dynamics transmits information
about its future from observations made in the past and the present.
Of particular interest are the long-term future, the periodicity and
so on. To capture salient features, there are essentially two
approaches: substantive and black-box. Examples of both
approaches abound. The former is often preferred if available in the
context in which we find ourselves. If not available, then a
black-box approach might be the only choice. We shall include
examples of both approaches. Let us first mention two substantive
examples as they are relevant to our later discussion.

\subsection{Two Substantive Models and Related Features}

1. \textit{Animal populations}. There are numerous ecological models
describing the time evolution of animal populations. The
single-species model of Oster and Ipaktchi (\citeyear{OstIpa78}) can be written as
\begin{equation}\label{animal}
\frac{dx_t}{dt} = b(x_{t-\tau}) x_{t-\tau} - \mu x_t,
\end{equation}
where $ x_t $ is the number of adults at time $ t $; $ \tau$ is the
delayed regulation duration due to the time taken for the young to
develop into adults or discrete breeding seasons; $b(\cdot)$ is the
birth rate; and $ \mu$ is the death rate. There are different
specifications for $b(\cdot) $. Gurney, Blythe and
Nisbet (\citeyear{GurBlyNis80}) suggested $
b(u) =\break c \exp(-u/N_0) $, where $ N_0 $ is the reciprocal of the
exponential decay rate and $ c$ is a parameter related to the
reproductive rate of adults. Ellner, Seifu and Smith (\citeyear{autokey22}) investigated
the estimation of model (\ref{animal}) by replacing $ b(x_{t-\tau})
x_{t-\tau} $ and $ \mu x_t $ with unknown functions $ B(x_{t-\tau})
$ and $ D(x_t) $, respectively, which they then used a nonparametric
method to estimate. Wood (\citeyear{Woo01}) considered a similar approach. There
are several discrete-time versions of (\ref{animal}) in biology.
See, for example, Varley, Gradwell and
Hassell (\citeyear{VarGraHas}). If we approximate $
dx_t/dt $ by $ x_t - x_{t-1} $, then\vadjust{\goodbreak} we obtain a nonlinear time
series model in discrete time
\begin{equation}\label{tpop1}
x_t = b(x_{t-\tau}) x_{t-\tau} + \nu x_{t-1}, 
\end{equation}
where 
$ \nu= 1-\mu$.

In ecology, population cycles are often observed and are an issue
of paramount importance. For example, the blowfly data show a cycle
of 39 days and the Canadian lynx shows a cycle of about 9.7 years.
Some ecologists have even suggested chaotic patterns, although we are
skeptical about this possibility. Most ecologists consider the
dynamics underlying population cycles as one of the major challenges
in their discipline.

2. \textit{Transmission of infectious diseases}. The conven\-tional
compartmental SIR model partitions a commu\-nity with population $ N $
into three compartments~$S_t$ (for susceptible), $I_t$ (for
infectious) and $R_t $ (for recovered): $ N = S_t + I_t + R_t $ at
any time instant $ t $. The SIR model is simple but very useful in
investigating many infectious diseases including measles, mumps,
rubella and SARS. Each member of the population typically
progresses from susceptible to infectious to recovered or death.

Infectious diseases tend to occur in cycles of outbreaks due to the
variation in the number of susceptible individuals over time. During
an epidemic, the number of susceptible individuals falls rapidly as
more of them are infected and thus enter the infectious and
recovered compartments. The disease cannot break out again until the
number of susceptible has built back up as a result of babies being
born into the susceptible compartment.

Consider a population characterized by a death rate~$ \mu$ and a
birth rate equal to the death rate, in which an infectious disease
is spreading. The differential equations of the SIR model are
\begin{eqnarray*}
\frac{dS_t}{dt} &=& \mu( N- S_t) - \beta\frac{I_t}{N} S_t,\\
\frac{dI_t}{dt} &=& \beta\frac{I_t}{N} S_t - (\nu+\mu) I_t,\\
\frac{dR_t}{dt} &=& \nu I_t - \mu R_t,
\end{eqnarray*}
where $ \beta$ is the contact rate and $\nu$ is the recovery rate of
the disease.
See, for example, Anderson and May (\citeyear{AndMay91}) and Isham and Medley (\citeyear{IshMed08})
for details.
This model has been extensively investigated and very successfully
used in the control of infectious diseases. Discrete-time versions
of the model have been proposed. An example is
\[ \label{tsir0}
I_{t+1} = r_0 S_t I_t/N,\quad  S_{t+1} = S_t - I_{t+1} + \mu N,
\]
where $ \mu$ is the birth rate and $ r_0 $ is the basic
reproductive rate of transmission. See, for example, Bartlett (\citeyear{BarN1}, \citeyear{BarN2}),
Anderson and May (\citeyear{AndMay91}) and the discussion in Section \ref{secDATA}.

Again, an important feature for the transmission of infectious
disease is the periodicity, to understand which it is essential to
understand the effect of such factors as the birth rate, the
seasonal force, the transmission rate and the incubation time on the
dynamics, the phase difference that is related to the transmission
in different areas, and the interaction between different diseases;
see, for example, Earn et al. (\citeyear{Earetal00}) and Rohani et al.
(\citeyear{Rohetal03}). The model can also be used to guide the policy maker in
controlling the spread of the disease. See, for example, Bartlett
(\citeyear{BarN1}), Hethcote (\citeyear{Het76}), Keeling and Grenfell (\citeyear{KeeGre97}) and Dye and
Gay (\citeyear{DyeGay}).

\vspace*{3pt}\subsection{The Objectives}\vspace*{3pt}

Our primary concern is parametric time series mo\-deling with the
objective of achieving good \textit{matching} of the joint probabilistic
distribution of the observable time series, including, in particular,
salient features, such
as cycles and others. Short-term
prediction is secondary in this paper. Accepting G.~E.~P. Box's  (\citeyear{Box76})
dictum: \textit{All models are wrong, but some are useful}, we
use parametric time series models only as means to an end.
We are typically less interested in the
consistency of estimators of unknown parameters in the conventional
sense, which is predicated on the assumed truth of the postulated
model. In fact, we are more interested in improving the
matching capability of
the postulated model.

Suppose we postulate the following model:
\begin{equation}\label{skeleton}
x_t = g_\theta(x_{t-1},\ldots, x_{t-p}) + \varepsilon_t,
\end{equation}
where $ \varepsilon_t $ is the innovation
and the function $ g_\theta(\cdot)$ is known up to parameters $ \theta$.
To indicate
the dependence of $ x_t $ on $ \theta$, we also write it as $
x_t(\theta) $. Following Tong (\citeyear{Ton90}), we call (\ref{skeleton}) with
$ \Var(\varepsilon_t ) = 0 $ the \textit{skeleton} of the model. In
postulating the above model, we recognize that it is generally just
an approximation of the true underlying dynamics no matter how the
function $ g_\theta(\cdot) $ is specified. Of particular note is the fact
that conventional methods of
estimation of~$\theta$ in the present setup are usually not
different
from those used for a cross-sectional model:
with observations $\{y_1, y_2, \ldots,
y_{_T}\}$ and postulated model $g_\theta$, typically a~loss function
is based on the
errors and takes the following
form:
%
\[
L(\theta) = (T-p)^{-1} \sum_{t=p+1}^T \{ y_t - g_\theta(y_{t-1},\ldots, y_{t-p})
\}^2,
\]
where, here and elsewhere, $T$ denotes the sample size. The errors above
happen to coincide with the one-step-ahead prediction errors. Under
general conditions, minimizing this loss function is known
mathematically to
lead to efficient estimation \textit{if the postulated model is true}.
However, the postulated model is, by the Box dictum, almost
invariably wrong, in which case the above loss function is not
necessarily fit for purpose.
To illustrate, let observations $\{ y_1, y_2,\ldots, y_T \}$ be given
and, of the postulated mo\-del~(\ref{skeleton}),
let the function $g_\theta$
be linear
and $\varepsilon_t$ be Gaussian with zero mean and finite
variance. Let $\mathcal{T} = \{C(j), j=0, 1, 2,\ldots, T-1\}$ denote a~set of
sample
autocovariances of the $y$-data. Then minimizing $L(\theta)$
yields well-known estimates of $\theta$ that are functions of $\mathcal{S}
=\{C(0), C(1),\ldots, C(p) \}.$ If the postulated model is ``right,''
then $\mathcal{S}$ is a minimal set of sufficient statistics (ignoring
boundary effects) and all is well. However, if it is wrong, then it
is unlikely that $\mathcal{S}$ will remain so. Since the model is typically
wrong, then restricting to $\mathcal{S}$ is unfit for the purpose of
estimating $\theta$; $\mathcal{T}$ may be preferable.

To reconcile with the Box spirit, diagnostic checks, goodness-of-fit
tests and other post-modeling devices are recommended. Indeed Box
and Jenkins (\citeyear{BoxJen70}) have stressed these post-modeling devices. See
also
Tsay (\citeyear{Tsa}) for some later developments. These are undoubtedly very
important developments. However, the challenge remains as to whether
we can adopt the Box spirit more seriously right at the modeling
stage rather than at the post-modeling stage.

It is worth recalling the fact that the classic autoregressive (AR)
model of Yule (\citeyear{Yul27}) and the moving average (MA) model of Slutsky
(\citeyear{Slu27}) were originally proposed to capture the sunspot cycle and the
business cycle, respectively, rather than for the purpose of short-term
prediction.

\section{The Matching Approach}

We shall use the letters $y$ and $x$ to signify respectively the real
time series under study
and the time series generated by the postulated model. The adjective observable
is reserved for a stochastic process. An observed time series
consisting of observations constitutes (possibly part of) a realization
of a stochastic process.
In order for model (\ref{skeleton}) to be able to
approximate an observable $ \{y_t\dvtx t = 1, 2,\ldots\} $ well, it is
natural to require throughout this paper that the state space of $
\{x_t(\theta)\dvtx t = 1, 2,\ldots\}$ covers that of the observable $
\{y_t\dvtx
t = 1, 2,
\ldots\} $. For expositional simplicity, let $ p = 1 $. Starting from
$ x_{0}(\theta) = y_{0}$, the postulated model is said to match an
observable
time series under study perfectly if their conditional distributions
are the same, namely,
\begin{eqnarray}\label{same}
&&\hspace*{25pt}P\{x_1(\theta_0)<u_1, \ldots, x_n(\theta_0)< u_n | x_{0}(\theta_0) =
y_0 \}\nonumber\\ [-8pt]\\ [-8pt]
&&\hspace*{25pt}\quad \equiv P\{y_1<u_1,\ldots, y_n<u_n |y_0\}\nonumber
\end{eqnarray}
almost surely for some $ \theta_0 $ and any $n$ and any real values
$ u_1,\ldots, u_n$. We call the approach based on the above model,
including all its weaker versions, some of which will be described
in the next two subsections, collectively by the name
\textit{catch-all approach}.

However, formulation (\ref{same}) is usually quite \mbox{difficult} to
implement in
practice. In the next two subsections, we suggest two weaker forms,
although other forms are obviously possible.

In the econometric literature, the notion of calibra\-tion has been introduced
(e.g., Kydland and Prescott, \citeyear{KydPre96}). It has many alternative
definitions.
Broadly speaking, calibration consists of a series of steps intended to
provide quantitative answers to a
particu\-lar economic question. A crucial step involves some so-called
``computational experiments''
with a substantive model of relevance to economic theory; it is
acknowledged that the model is
unlikely to be the true model for the observed economic data. At the
philosophical level,
calibration and our feature matching share almost the same aim.
However, there are some fundamental differences
in methodology. Our methodology provides a~statistical and \textit{coherent} framework (in a~non-Bayesian sense) to estimate \textit{all} the
parameters of a postulated (and usually wrong) model. As far as we
know, calibration seems to be in need of such a framework.
See, for example, Canova (\citeyear{Can07}), esp. page 239. The hope is that our
methodology will be useful to
substantive modelers in all fields, including ecology, economics,
epidemiology and others. At the other end
of the scale, it has been suggested that our methodology has potential
in data mining
(K.S. Chan, private communication).

\subsection{Matching Up-to-$m$-Step-Ahead Point Predictions}\vspace*{3pt}

If we are interested in the mean conditional on some
initial observation, say $y_0$,
we can weaken the matching
requirement (\ref{same}) to
\begin{eqnarray*} \label{expectation}
&&E[(x_1(\theta_0),\ldots, x_m(\theta_0))|x_{0}(\theta_0) =
y_{0}]\\
&&\quad \equiv E[(y_1,\ldots, y_m)| y_{0}],
\end{eqnarray*}
where the length $m$ of the random vector is, in practice, bounded
above by the sample size
under consideration. The expectation is taken with respect
to the relevant joint distribution of the random vector conditional on
the initial
value being $y_0.$
Since a postulated model is just an approximation of the underlying
dynamics, we set $ \theta_{0} $ to minimize the difference of the
prediction vectors, that is,
\begin{eqnarray}\label{expectA}
&&\hspace*{20pt}E\{\|E[(x_1(\theta),\ldots, x_m(\theta))|x_{0}(\theta)
= y_{0}]\nonumber\\ [-8pt]\\ [-8pt]
&&\hspace*{90pt}{} - E[(y_1, \ldots, y_m)| y_{0}]\|^2\}.\nonumber
\end{eqnarray}
Here, $\| \cdot\|$ denotes the Euclidean norm of a vector. In other
words, we choose $ \theta$ by minimizing up-to-$m$-step-ahead
prediction errors. It is basically based on a \textit{catch-all} idea. It
is easy to see that the best $ \theta$ based
on minimizing (\ref{expectA}) depends on $ m $.
Generally speaking, we set $m = 1,$ when and only when we have
complete faith in the model, which is what the conventional methods
do. Denote the $m$-step-ahead prediction of $ y_{t+m} $ based on
model (\ref{skeleton}) by
\[
g_\theta^{[m]}(y_t) = \E( x_{t+m}|x_t = y_t).
\]
If model (\ref{skeleton}) is
deterministic [i.e., $ \Var(\varepsilon_t) =0 $] or linear,
$g_\theta^{[m]}(y_t ) $ is simply a composite function,
\[
g_\theta^{[m]}(y_t) = \underbrace{g_\theta(g_\theta(\cdots g_\theta(}_{m\ \mathrm{folds}}y_t)\cdots)).
\]

Let
\begin{eqnarray}\label{criterion}
&&Q(y_t, x_t(\theta))\nonumber\\ [-8pt]\\ [-8pt]
&&\quad= \sup_{w_m} \sum_{m=1}^{\infty} w_m \bigl[ \E
\bigl\{ y_{t+m} - g_\theta^{[m]}(y_t)\bigr\}^2\bigr],\nonumber
\end{eqnarray}
where $ w_m \ge0 $ and $ \sum w_m = 1 $. Since
\[
\E\bigl[ \{ y_{t+m} - \E(y_{t+m}|y_t) \} \bigl\{\E(y_{t+m}|y_t)- g_\theta
^{[m]}(y_t)\bigr\}\bigr] = 0,
\]
we have
\begin{eqnarray*}
\E\bigl\{ y_{t+m} - g_\theta^{[m]}(y_t)\bigr\}^2 &=& \E\{ y_{t+m} - \E(y_{t+m}|y_t) \}^2\\
&&{} + \E\bigl\{\E(y_{t+m}|y_t)- g_\theta^{[m]}(y_t)\bigr\}^2.
\end{eqnarray*}
Let
\begin{eqnarray*}
&&\tilde Q(y_t, x_t(\theta))\\
&&\quad = \sup_{w_m}\sum_{m=1}^\infty w_m \bigl[
\E\bigl\{\E(y_{t+m}|y_t)- g_\theta^{[m]}(y_t)\bigr\}^2\bigr].
\end{eqnarray*}
If the observable $ y_t $ indeed follows the model of~$ x_t $, then
$\min_\theta\tilde Q(y_t, x_t(\theta)) = 0 $. Otherwise we generally
expect $\min_\theta\tilde Q(y_t, x_t(\theta))
> 0 $. Minimizing $ \tilde Q(y_t, x_t(\theta)) $ is for $x_t(\theta) $ to arrive
at a choice within the postulated model that
gives all (suitably weighted) multiple-step-ahead predictions of $
y_t $ as accurately as possible in the mean squared sense.

Note that the above measure of the difference between two time
series is based on a (weighted) least squares loss function. Clearly
there exist many other possible measures. For example, if the
distribution of the innovation is known, a likelihood type measure
of the difference can be used instead. A Bayesian may perhaps then
endow $\{w_m\}$ with some prior distribution. This line of
development may be worth further exploration as suggested by an
anonymous~re\-feree. Intuitively speaking, a $J$-shaped $\{w_m\}$
tends~to emphasize low-pass filtering, because $\E(y_{t+m}|y_t)$ is
a~slowly varying function for sufficiently large $m$.~Si\-milarly, an
inverted-$J$-shaped $\{w_m\}$ tends to empha\-size high-pass
filtering. An optimal choice of $\{w_m\}$ strikes a good balance
between high-pass filtering and low-pass filtering.

The most commonly used estimation method in time series modeling is
probably that based on minimizing the sum of squares of the errors
of one-step-ahead prediction.
This has been extended to the sum of squares of errors of other
single-step-ahead prediction. See, for example, Cox (\citeyear{Cox61}), Tiao and Xu
(\citeyear{TiaXu93}), Bhansali and Kokoszka (\citeyear{BhaKok02}) and Chen, Yang and
Hafner (\citeyear{CheYanHaf04}). Clearly,
the former method is predicated on the
model being \textit{true}. The latter extension recognizes that this is an
unrealistic assumption for multi-step-ahead prediction. Instead, a
\textit{panel of models} is constructed so that a different model is
used for the prediction at each different horizon. The focus of
the extension is prediction.

The approach that we develop here
essentially\break builds on the above extension. First, we shift the
focus away from prediction. Second, we
transform the prediction based on
a \textit{panel of models} into the fitting of a \textit{single time series model}.
We effectively synthesize the
panel into a catch-all methodology. Specifically, we propose to
minimize the sum of
squares of errors of prediction \textit{over all \textup{(}allowable\textup{)} steps
ahead}, as given in (\ref{criterion}). We
stress again that our primary objective is feature matching
rather than prediction. Of course, it is conceivable that good
feature matching may sometimes lead to better prediction, especially
for the medium and long term. Clearly each member of the panel can
be recovered, at least formally, from the catch-all setup by
setting, in turn, the weight, $w_j$, to unity, leaving the rest to
zero.\

\subsection{Matching ACFs}\label{secKK}

Suppose that the observable $ \{y_t \} $ and $ \{ x_t(\theta) \} $ are both
second-order stationary. If we are interested in second-order
moments, then a weaker form of (\ref{same}) is the
following difference or distance function:
\[
D_{_C}(y_t, x_t(\theta)) = \sup_{\{w_m\}}\sum_{m=0}^\infty w_m \bigl\{
\gamma_{x(\theta)}(m) - \gamma_y(m) \bigr\}^2.
\]
Here, the suffixes of
$y$ and
$x(\theta)$ are self-explanatory. We assume that the spectral
density function (SDF) of the observable $ y_t $ exists; it is given by
\[
f_{y}(\omega) = \frac1{2 \pi}\gamma(0) + \frac1\pi\sum
_{k=1}^\infty\gamma_{y}(k) \cos(k
\omega).
\]
The SDF of $x_t (\theta)$, which we also assume to exist, can be
defined similarly. We can also measure the difference between two
time series by reference to the difference between their SDFs, for
example,
\[
D_{_F}(y_t, x_t(\theta)) = \int_{-\pi}^{\pi}\biggl\{\frac{f_{y}(\omega
)}{f_{x}(\omega)} + \log\biggl(\frac{f_{x}(\omega)}{f_{y}(\omega)}\biggr) - 1
\biggr\}\,d\omega,
\]
which is called the Itakura--Saito distortion measure; see also
Whittle (\citeyear{Whi62}). Further discussion on measuring the difference
between two SDFs can be found in Georgiou (\citeyear{Geo07}).


Suppose that $\{x_t(\theta)\}$ and the observable $\{y_t\}$ have the
same marginal
distribution and they each have second-order moments. Then we can prove
that
\begin{eqnarray*}
D_{_C}(y_t, x_t(\theta)) &\le& C_1 \tilde Q(y_t, x_t(\theta)),\\
D_{_F}(y_t, x_t(\theta)) &\le& C_2 \tilde Q(y_t, x_t(\theta))
\end{eqnarray*}
for some positive constants $ C_1 $ and $C_2$. Moreover, if
$\{x_t(\theta)\}$ and the observable $\{y_t\}$ are linear AR models,
then there are some positive constants $ C_3 $ and~$ C_4 $ such that
%
\begin{eqnarray*}
\tilde Q(y_t, x_t(\theta)) &\le& C_3 D_{_C}(y_t, x_t(\theta)),\\
\tilde Q(y_t, x_t(\theta)) &\le& C_4 D_{_F}(y_t, x_t(\theta)).
\end{eqnarray*}
For further details, see Theorem \ref{thma} in the \hyperref[appendix]{Appendix}.

For linear AR models under the above setup,
$ \tilde Q(\cdot,\cdot)$, $ D_C(\cdot, \cdot) $ and $D_F(\cdot,\cdot)$ are equivalent.
However, the equivalence is not generally true.
A counterexample can be constructed easily by reference to the
classic random telegraph signal process. [See, e.g., Parzen (\citeyear{Par62}),
page 115].

Let us close this section by describing one way of implementing the
ACF criterion
for an ARMA model with
normal innovation. Suppose $ y_1,\ldots, y_{_T} $ are observations from
the observable $
\{ y_t \} $. Whittle (\citeyear{Whi62}) considered a ``likelihood function'' for
ARMA models in terms of the SDF. Let
\[
I(w) =\frac{1}{2\pi T} \Biggl|\sum_{t=1}^T y_t \exp(- \iota\omega t)\Biggr|^2
\]
be the periodogram of the sample, where $ \iota$ is the imaginary
unit. Let $ f_\theta(\omega) $ be the theoretical SDF of an ARMA model
with parameters $ \theta$. Whittle (\citeyear{Whi62}) proposed to estimate $
\theta$ by
\[
\hat\theta= \min_{\theta} \sum_{j=1}^T \biggl\{ \frac{I(\omega
_j)}{f_\theta(\omega_j)} + \log(f_\theta(\omega_j)) \biggr\},
\]
where $ \omega_j = 2\pi j/T $.
From the
perspective of feature matching, the celebrated Whittle's likelihood
is not a conventional likelihood
but a precursor of the ex\-tended-likelihood approach.
It matches
the second-order moments, by using a natural sample version of
$D_F(y_t, x_t(\theta))$ up to a constant.
For this reason, it is expected that for misspecified models,
Whittle's estimator can lead to better matching of the ACFs of the
observed time series than the innovation driven methods [e.g., the
least squares estimation (LSE) or the maximum likelihood estimation
(MLE)]. We shall give some numerical comparison between Whittle's
estimator and the others in Sections \ref{sec5} and \ref{secDATA} below.


\section{Time Series with Measurement Errors}

To illustrate the advantages of the catch-all approach, which
involves minimal assumptions on the observed time
series, we give detailed analyses of two cases involving
measurement errors, one of which
is related to a linear
$\operatorname{AR}(p)$ model and the other a~nonlinear skeleton model. They
can be considered special cases of model misspecification in that
the observable $y$-time series is a measured version of the $x$-time
series subject
to measurement errors.
For the linear case, measurement error is an old
problem in time series analysis that was studied at least as early as
Walker (\citeyear{Wal60}). Some
new lights will be shed.

\subsection{Linear $\operatorname{AR}(p)$ Models} \label{SecL}

Consider the following $\operatorname{AR}(p)$ model:
%
\begin{equation} \label{ssm0}
x_{t} = \theta_1 x_{t-1} + \cdots+ \theta_p x_{t-p} + \varepsilon_{t}.
\end{equation}
Stationarity is assumed.
By the Yule--Walker
equations,
we have the recursive formula for the
ACF, $\{\gamma(j)\}$, of the $x$-time series, namely,
\begin{eqnarray}\label{ywk}
\gamma(k) &=& \gamma(k-1)\theta_1 + \gamma(k-2)\theta_2 + \cdots\nonumber\\ [-8pt]\\ [-8pt]
&&{}+
\gamma(k-p)\theta_p,\quad
k = 1, 2, \ldots.\nonumber
\end{eqnarray}
Let $m \ge p$ and $\Upsilon_m = (\gamma(1), \gamma(2),\ldots,
\gamma(m))^\top$, $ \theta= (\theta_1,\ldots, \theta_p)^\top$ and
\[
\Gamma_m =
\pmatrix{\gamma(0) & \gamma(-1) & \cdots &\gamma(-p+1) \cr
\gamma(1) & \gamma(0) & \cdots &\gamma(-p+2) \cr
\vdots\cr
\gamma(m-1) & \gamma(m-2) & \cdots &\gamma(m-p)
}
.
\]
The Yule--Walker equations can be written as
\[
\Gamma_m \theta=
\Upsilon_m .
\]

Suppose that the observable $y$-time series is given by $y_t = x_t +
\eta_t$,
for $t = 0, 1, 2, \ldots,$ where $\{\eta_t\}$ is independent of $\{
x_t\}$ and is
a sequence of independent
and identically distributed random variables each with zero mean and finite
variance. Clearly, $\{y_t\}$ is no longer given by an $\operatorname{AR}(p)$ model of the
form (\ref{ssm0}).

Let $\{\tilde\gamma(j) \}$ denote the ACF of the observable
$y$-time series. Let $ \tilde\Gamma_m $ and $ \tilde\Upsilon_m $
denote the analogously defined matrix and vector of ACFs for the
observable $y$-time series.

Suppose we are now given the observations
$\{y_1, y_2,\allowbreak \ldots, y_T\},$ and we wish to
fit the wrong model of the form (\ref{ssm0}) to them.
We may estimate $\tilde\gamma(j)$ by
$ \hat\gamma(j) =
\hat\gamma(-j) = T^{-1}\sum_{t=1}^{T-j} (y_t -\bar y)(y_{t+j} -\bar
y) $, $\bar y$ being the sample mean. Let
$\hat\Gamma_m $ and $ \hat\Upsilon_m $
denote the obvious sample version of $\tilde\Gamma_m $
and sample version of $\tilde\Upsilon_m $, respectively.


Since any $ p $ equations can be used to determine
the parameters, the Yule--Walker estimators typically use the first
$p$ equations, that is,
\[
\hat\theta= \hat\Gamma_p^{-1} \hat\Upsilon_p \quad\mbox{or}\quad\hat\theta= (\hat\Gamma_p^\top\hat\Gamma_p)^{-1} \hat
\Gamma_p^\top\hat\Upsilon_p,
\]
which is also the minimizer of $ \sum_{k=1}^p \{ \hat\gamma(k) -
\hat\gamma(k-1)\theta_1 - \hat\gamma(k-2)\theta_2 - \cdots-\hat
\gamma(k-p)\theta_p \}^2 $, involving the ACF only up to lag $p$. We
can achieve closer matching of the ACF by incorporating lags beyond
$p$ as well. For example, we may consider estimating\vadjust{\goodbreak} $ \theta$ by
minimizing
\begin{eqnarray}
&&\sum_{k=1}^{m} \{ \hat\gamma(k) - \hat\gamma(k-1)\theta_1\nonumber\\
&&\hphantom{\sum_{k=1}^{m} \{}{} - \hat \gamma(k-2)\theta_2 - \cdots- \hat\gamma(k-p)\theta_p \}^2,\qquad
\eqntext{m \ge p.}
\end{eqnarray}
Denoting the minimizer by $ \hat\theta_{\{m\}} $, we have
\begin{eqnarray}
\hspace*{15pt}\hat\theta_{\{m\}} &=& (\hat\Gamma_{m}^\top\hat\Gamma_{m})^{-1}
\hat\Gamma_{m}^\top\hat\Upsilon_{m}\nonumber\\ [-8pt]\\ [-8pt]
\hspace*{15pt}&=&
\Biggl\{\sum_{k=0}^m \breve\Upsilon_k \breve\Upsilon_k^\top\Biggr\}^{-1}
\sum_{k=0}^m \breve\Upsilon_k \hat\gamma(k+1),\nonumber
\end{eqnarray}
where $\breve\Upsilon_k = (\hat\gamma(k), \hat\gamma(k+1),\ldots
, \hat\gamma(k+p-1))^\top.$
Let us call the estimator $ \hat\theta_{\{m\}} $ the
up-to-lag-$m$ Yule--Wal\-ker estimator (or AYW($\le m)$).
For the error-free case, that is,
$\eta_t = 0 $ with probability 1, it is easy to see that $ \hat\theta
_{\{p\}} $ is
the most efficient amongst all~$ \hat\theta_{\{m\}}$, $m = p, p+1,
\ldots. $ Otherwise,
under some regularity conditions,
we have in
distribution
\[
\sqrt{n} \bigl\{ \hat\theta_{\{m\}} - \vartheta\bigr\} \to N(0, \tilde\Sigma_m),
\]
where $ \vartheta= (\tilde\Gamma_{m}^\top\tilde\Gamma_{m})^{-1}
\tilde\Gamma_{m}^\top\tilde\Upsilon_{m}$ and $ \tilde\Sigma_m $
is a positive~de\-finite matrix.
For $ \Var(\varepsilon_t) >0$ and \mbox{$\Var(\eta_t) = \sigma_\eta^2>0 $}, the above asymptotic result holds with
$ \vartheta= \theta+\break\sigma_\eta^2 (\Gamma_m^\top\Gamma_m + 2
\sigma_\eta^2\Gamma_p + \sigma_\eta^4 I)^{-1} (\Gamma_p +
\sigma_\eta^2 I)\theta. $
For further details, see
Theorem \ref{thmb} in the \hyperref[appendix]{Appendix}.


Clearly the bias $\sigma_\eta^2 (\Gamma_m^\top\Gamma_m + 2 \sigma
_\eta^2\Gamma_p + \sigma_\eta^4 I)^{-1} (\Gamma_p + \sigma_\eta
^2 I)\theta$ in the estimator will be smaller when $ m $ is~lar\-ger. For
sufficiently large sample size, the smaller~bias can lead to higher
efficiency in the sense of~mean squared errors (MSE). Let $ \bar
\Upsilon_k = (\gamma(k), \gamma(k+1),\ldots,\allowbreak \gamma(k+p-1))^\top$.
Then
\[
\Gamma_m^\top\Gamma_m = \Gamma_p^\top\Gamma_p + \sum_{k=p}^m
\bar\Upsilon_k \bar\Upsilon_k^\top.
\]
Thus, the bias can be reduced more substantially if the ACF decays very
slowly and a larger $ m $ is used. For example, a highly cyclical
time series usually has slowly decaying ACF, in which case the AYW
will provide a substantial improvement over the Yule--Walker
estimators. However, even with the ACF slowly decaying, a large $ m$
may cause larger variability of the estimator. Therefore, a good
choice of $ m $ is also important in practice. We shall return to
this issue later.

In fact, Walker (\citeyear{Wal60}) suggested using exactly $ p $ equations to
estimate the coefficients giving
\[
\hat\theta_{W.\ell} = \operatorname{arg\,min}\limits_\theta\Biggl\{\sum_{k=p+\ell
}^{2p-1+\ell} \breve\Upsilon_k \breve\Upsilon_k^\top\Biggr\}^{-1}
\sum_{k=p+\ell}^{2p-1+\ell} \breve\Upsilon_k \hat\gamma(k+1).
\]
Note the difference between AYW and $ \hat\theta_{W.\ell} $. Walker
(\citeyear{Wal60}) showed that in the presence of measurement error, then $ \ell
= p $ is the optimal choice amongst all candidates with $\ell\ge p
$, by reference to MSE. However, Walker's method seems
counterintuitive because it relies on the sample ACF at higher lags
to a greater extent than those at the lower lags.
Further discussion on Walker's method can be found in Sakai, Soeda and
Tokumaru (\citeyear{SakSoeTok79}) and Staudenmayer and Buonaccorsi (\citeyear{StaBuo05}). It is
well known that an autoregressive model plus independent additive
white noise results in an ARMA model. Walker's approach essentially
treats the resulting ARMA model as a true model. This approach has
attracted attention in the engineering literature. See, for example,
Friedlander and Sharman (\citeyear{FriSha85}) and Stoica, Moses and Li (\citeyear{Sto}).
The essential difference between this approach and the catch-all
approach is that the latter postulates an autoregressive model to
match the observations. And we know that it is a wrong model, as we
consistently do with all postulated models. Note that the use of
sample ACFs at all possible lags has points of contact with the
so-called generalized method of moments, used extensively in
econometrics. See, for example, Hall (\citeyear{Hal05}).

Next, we consider estimation based on $Q(\cdot,\cdot) $. Gi\-ven a finite
sample size, we may stop at, say, the $m$-step-ahead prediction.
Let $ e_1 = (1, 0,\ldots, 0)^\top
$ and
\[
\Phi=
\pmatrix{ \theta_1 & \theta_2 & \cdots & \theta_{p-1} & \theta_p
\cr
1 & 0 & \cdots & 0 & 0 \cr
\vdots& \vdots& \vdots& 0 & 0 \cr
0 & 0 & \cdots & 1 & 0
}.
\]
We estimate $ \theta$ by
\begin{eqnarray}\label{mL}
\tilde\theta_{\{m\}}= \operatorname{arg\,min}\limits_{\theta} \sum_{t=p+1}^T
\sum_{k=1}^m w_k \{ y_{t-1+k}\nonumber\\ [-8pt]\\ [-8pt]
\eqntext{ - e_1^\top\Phi^k (y_{t-1},\ldots,
y_{t-p})^\top\}^2,}
\end{eqnarray}
where $ w_k $ is a weight function, typically positive definite. A
reasonable choice of $ w_k$ is the absolute value of the
autocorrelation function of the observed time series,
that is, $ w_k = |r_y(k)|$. We call $ \tilde\theta_{\{m\}} $ in
(\ref{mL}) the up-to-$m$-step-ahead prediction estimator [$\operatorname{APE}$ or
$\operatorname{APE}(\le m)$].

The asymptotic properties of $ \tilde\theta_{\{m\}} $ will be
discussed later.


\subsection{Nonlinear Skeletons} \label{SecN}

A deterministic nonlinear dynamic model with measurement error is
commonly used in many applied areas, for example, ecology, dynamical
systems\vadjust{\goodbreak}
and others. See, for example, May (\citeyear{May76}),
Gurney, Blythe and Nisbet (\citeyear{GurBlyNis80}), Tong (\citeyear{Ton90}), Anderson and May (\citeyear{AndMay91}),
Alligood, Sauer and
Yorke (\citeyear{AllSauYor97}), Grenfell, Bj{\o}rnstad and
Finkenst{\"a}dt (\citeyear{GreBjrFin02}), Chan and
Tong (\citeyear{ChaTon01}) and the examples in Section \ref{secDATA}.
Consider using the following nonlinear skeleton:
\begin{equation}
x_t = g_\theta(x_{t-1},\ldots, x_{t-p}) \label{nonlinearmodel}
\end{equation}
to match the observable time series $\{y_t \} $.

Employing the $Q(\cdot,\cdot)$ criterion, the estimator is given by
\begin{eqnarray} \label{mN}
\tilde\theta_{\{m\}} = \operatorname{arg\,min}\limits_{\theta} \sum_{t=p+1}^T
\sum_{k=1}^m w_k \bigl\{ y_{t-1+k}\nonumber\\ [-8pt]\\ [-8pt]
\eqntext{ - g^{[m]}_\theta( y_{t-1},\ldots,
y_{t-p}) \bigr\}^2,}
\end{eqnarray}
which we again call the up-to-$m$-step-ahead prediction estimator
[$\operatorname{APE}$ or $\operatorname{APE}(\le m)$]. Here the weight function $\{w_k\}$ is as
defined in (\ref{criterion}).

For ease of explanation, we consider again $ y_t = x_t + \eta_t
$ and $ p = 1 $.
Starting from any state $ \tilde x_0 = x_0 $, let $ \tilde x_t =
g_\theta^{[m]}(x_0) $. Suppose the dynamical system has a negative
Lyapunov exponent
\[
\lambda_\theta(x_0) = \lim_{n\to\infty} n^{-1} \sum
_{t=0}^{n-1}\log( |g_\theta'(\tilde x_t)|) < 0,
\]
for all states $ x_0 $. Similarly starting from $ x_t $ let
$
x_{t+m} = g^{[m]}_{\theta_0}( x_t)
$.
We
predict $ x_{t+m} $ by $
\hat y_{t+m} = g^{[m]}_\theta(y_t)
$. By the definition of the Lyapunov exponent, we have
\[
\bigl| g^{[m]}_\theta(x_t+ \eta_t) - g^{[m]}_\theta(x_t)\bigr| \approx\exp\{
m \lambda_\theta(x_t)\} |\eta_t|.
\]
%
More generally, suppose the
system $ x_t=g_{\theta_0}(x_{t-1},\allowbreak\ldots, x_{t-p}) $ has a
finite-dimensional state space and admits only limit cycles,
but $x_t$ is observed
as $ y_t = x_t + \eta_t $, where $ \{\eta_t \}$ are
independent with mean 0. Suppose that the function $ g_\theta(v_1,
\ldots, v_p) $ has bounded derivatives in both $ \theta$ in the
parameter space $\Theta$ and $ v_1, \ldots, v_p $ in a neighborhood
of the state space. Suppose that the system $ z_t =
g_\theta(z_{t-1},\ldots, z_{t-p}) $ has only negative Lyapunov
exponents in a small neighborhood of $ \{x_t\} $ and in $ \theta
\in\Theta$. Let $ X_t = (x_t, x_{t-1},\ldots,\allowbreak x_{t-p}) $ and $ Y_t
= (y_t, y_{t-1},\ldots, y_{t-p}) $.
If the observed
$Y_0 = X_0 + (\eta_{0}, \eta_{-1},\ldots, \eta_{-p}) $
is taken as the initial values of $ \{x_t\}$, then for any $n$,
\begin{eqnarray}\label{result001}
&&f(y_{m+1}, \ldots, y_{m+n}|X_0 )\nonumber\\
&&\quad{} - f(y_{m+1}|X_0 = Y_0)\\
&&\qquad{}\cdots
f(y_{m+n}|X_0 = Y_0) \to0\nonumber
\end{eqnarray}
as $ m \to\infty$. Suppose the equation $ \sum_{X_{t-1}} \{ g_\theta
(X_{t-1}) -x_t\}^2 = 0 $\vadjust{\goodbreak}
has a unique solution in $\theta$,
where the summation is taken over all limiting states. Let $
\theta_{\{m\}} =
\operatorname{arg\,min}_{\theta} m^{-1} \sum_{k=1}^m \EE\{ y_{t-1+k}
- g^{[k]}_\theta( Y_{t-1}) \}^2. $ If the\break noise takes value in
a small neighborhood of the origin, then
\[
\theta_{\{m\}} \to\theta_0\quad \mbox{as $ m \to\infty. $}
\]
%
Note that $ |f(y_{1}|X_0) -
f(y_{1}|X_0 = Y_0) | \neq0 $ implies that
\begin{eqnarray*}
&&f(y_{1},\ldots, y_n|X_0= Y_0)\\
&&\quad \neq
f(y_{1}|X_0=Y_0)f(y_{2}|X_1=Y_1)\\
&&\qquad\cdots
f(y_{n}|X_{n-1}=Y_{n-1}),
\end{eqnarray*}
which
challenges the commonly used (conditional) MLE. Equation (\ref
{result001}) indicates that using high step-ahead
prediction can reduce the effect of
noisy data (e.g., due to measurement errors), and provide
a better approximation of the conditional distribution. The second
part suggests that using high step-ahead prediction errors in a
criterion can reduce the bias caused by the presence of
$\eta_t$. It
also implies that any set of past values, for example, $ (y_{t-1},\ldots,
y_{t-p}) $ for $ t> p $, can offer us an estimator with the first
summation in (\ref{mN}) removed. However, the summation over all
past values is more efficient statistically.
For further details, see Theorem \ref{thmc} in the \hyperref[appendix]{Appendix}.

There are other interesting special cases. For example, when the
postulated model has a chaotic skeleton, the initial values play a
crucial role. One approach is to treat the initial values as unknown
parameters. See, for example, Chan and Tong (\citeyear{ChaTon01}) for more details.
Another example is when the postulated model is nonlinear, and is
driven by nonadditive white noise with an unknown distribution.
Here, the exact least squares multi-step-ahead prediction is quite
difficult to obtain theoretically and time consuming to calculate
numerically; see, for example, Guo, Bai and An (\citeyear{GuoBaiAn99}). In this case, the
up-to-$m$-step-ahead prediction method is difficult to implement
directly. However, our simulations suggest that approximating the
multi-step-ahead prediction by its skeleton is sometimes helpful in
feature matching, especially when the observed time series is quite cyclical
(Chan, Tong and Stenseth, \citeyear{ChaTonSte09}).


\section{Issues of the Estimation Method}

We now turn to some theoretical issues and calculation problems. In
conventional statistical theory for parameter estimation, by
consistency is generally meant that the estimated parameter vector
converges to the true parameter vector in some sense as the sample
size tends to infinity. The postulated model is assumed to be
the true model in the above conventional approach.

In the absence of a true model and \textit{ipso
facto} true parameter vector, we propose an alternative
definition of consistency. Specifically, by consistency we mean
that the estimated parameter vector will, in some sense, tend to the
optimal parameter vector that represents the best achievable
feature matching
of the postulated model to the observable time series.
To be more precise, for some positive integer~$m$
(which may be infinite), we define the optimal parameter by
\begin{eqnarray*}
&&\vartheta_{m, \mathbf{w}} = \operatorname{arg\,min}\limits_\theta\sum_{k=1}^m w_k \EE
[ y_{t+k}\\
&&\hphantom{\vartheta_{m, \mathbf{w}} = \operatorname{arg\,min}\limits_\theta\sum_{k=1}^m w_k \EE} {}-
\EE\{x_{t+k}(\theta)|X_t(\theta)=Y_t\}]^2,
\end{eqnarray*}
where $ X_t(\theta) = (x_t(\theta),\ldots, x_{t-p+1}(\theta))$ and
$\{w_k\}$ defines the weight function, typically positive and
summing to unity. For ease of exposition, we assume that the
solution to the above minimization is unique. Now, we say that an
estimator is \textit{feature-consistent} if it converges to
$\vartheta_{m,\mathbf{w}}$ in probability as the sample size tends to
infinity. It is easy to prove that under
some regularity conditions, $ \tilde\theta_{\{m\}} $ is
asymptotically normal, that is,
\[
T^{-1/2} \bigl(\tilde
\theta_{\{m\}} - \vartheta_{m, \mathbf{w}} \bigr) \stackrel{D}{\to} N(0,
\Omega)
\]
for some positive definite matrix $ \Omega$.
For further details, see Theorem \ref{thmd} in the \hyperref[appendix]{Appendix}.

The optimal parameter depends on $ m $ and the weight function $
w_k $. As discussed in Section \ref{SecL}, when the autocorrelation
decays less slowly, we should consider using a larger $m$.
Alternatively, we can consider assigning heavier weights for larger
$k$. Our experience suggests that, for a postulated linear time
series model, $ w_k $ can be selected as the absolute value of the
sample ACF function. For a postulated nonlinear time series model
aiming to match possibly high degrees of periodicity, $ w_k $ can be
chosen as constant lasting for approximately one, two or three periods. Note
that by setting $w_1 = 1$ and all other~$w_j$'s zero, the estimation
is equivalent to the LSE, and the MLE in the case of exponential
family of distributions.

The above feature suggests that we may regard~$\tilde\theta_{\!\{\!m\!\}}$
as a \textit{maximum extended-likelihood estimator} and fun\-ctions such as
$\sum_{t=p+1}^T \sum_{k=1}^m w_k \{ y_{t-1+k} - e_1^\top\Phi^k
(y_{t-1},\allowbreak\ldots, y_{t-p})^\top\}^2$
or their equivalents as \textit{ex\-tended-likeli\-hoods} (or \textit{XT-likelihoods} for short), with Whittle's likelihood
as a precursor. An XT-likelihood carries with it the interpretation as
a weighted average of likelihoods of
a cluster of models around the postulated model. In this sense, it is
related to
Akaike's notion of the likelihood of a model (Akaike, \citeyear{Aka78}).

For the numerical calculation involved in (\ref{mL}) and (\ref{mN}),
the gradient and the Hessian matrix of the loss function can be
obtained recursively for different steps of prediction. Consider
(\ref{mN}) as an example. Let $ g_\theta^{[m]} $ stand for $
g_\theta^{[m]}(y_{t-1},\ldots, y_{t-p}) $ and write $ g_\theta(v_1,
\ldots, v_p) $ as $ g(v_1,\ldots, v_p, \theta_1,\ldots, \theta_q) $. Let $
g_\theta^{[0]} =\break y_{t-1},\ldots, g_\theta^{[-p+1]} = y_{t-p} $, $
{\partial g_\theta^{[m]}}/{\partial\theta_k} = 0 $ and $
{\partial^2 g_\theta^{[m]}}/\allowbreak({\partial\theta_k \partial
\theta_\ell})=0, k,\ell=1,\ldots, q $ if $ m \le0 $. Then for $ m
\ge1 $,
\[
g_\theta^{[m]} = g\bigl(g_\theta^{[m-1]},\ldots, g_\theta^{[m-p]}, \theta
_1,\ldots, \theta_q\bigr)
\]
and
\begin{eqnarray}
\frac{\partial g_\theta^{[m]}}{\partial\theta_k} &=& \sum_{i=1}^p
\dot{g}_{i} \frac{\partial g_\theta^{[m-i]}}{\partial\theta_k}\nonumber\\
&&{} +
\dot{g}_{p+k}\bigl(g_\theta^{[m-1]},\ldots, g_\theta^{[m-p]}, \theta_1,
\ldots, \theta_q\bigr),\nonumber\\
\eqntext{ k = 1,\ldots,q,}
\end{eqnarray}
where $ \dot{g}_k(v_1, \ldots, v_p, \ldots, v_{p+q}) = \partial g(v_1, \ldots,
v_p, \ldots,\break v_{p+q})/\partial v_k, k = 1, \ldots, p+q $, and
\begin{eqnarray*}
&&\frac{\partial^2 g_\theta^{[m]}}{\partial\theta_k \,\partial\theta
_\ell}\\
&&\quad= \sum_{i=1}^p \sum_{j=1}^p \ddot{g}_{i,j} \frac{\partial
g_\theta^{[m-i]}}{\partial\theta_k} \frac{\partial g_\theta
^{[m-j]}}{\partial\theta_\ell} + \sum_{i=1}^p \dot{g}_{i} \frac
{\partial^2 g_\theta^{[m-i]}}{\partial\theta_k\, \partial\theta
_\ell} \\
&&\qquad{}+ \sum_{i=1}^p \ddot{g}_{p+k,i}\bigl(g_\theta^{[m-1]}, \ldots, g_\theta
^{[m-p]},\\
&&\hspace*{127.2pt} \theta_1, \ldots, \theta_q\bigr)\frac{\partial g_\theta
^{[m-i]}}{\partial\theta_\ell}\\
&&\qquad{} +
\ddot{g}_{p+k, p
+\ell}\bigl(g_\theta^{[m-1]}, \ldots, g_\theta^{[m-p]}, \theta_1, \ldots,
\theta_q\bigr),
\end{eqnarray*}
where $ \ddot{g}_{k,\ell}(v_1, \ldots, v_p, \ldots, v_{p+q}) = \partial^2
g(v_1, \ldots, v_p, \ldots,\allowbreak v_{p+q})/(\partial v_k \partial v_\ell) $ for $
k, \ell= 1, \ldots,p+q $. The Newton--Raphson method can then be used
for the minimization.

\section{Simulation Study}\label{sec5}

There are many
different ways to measure the goodness of matching the observed by
the postulated model, depending on the features of interest. We
suggest two here. (1) The ACFs are clearly important features in the
context of linear time series, and relevant even for nonlinear time
series analysis. Therefore, a natural measure can be based on the
differences of the ACFs, for example,
\begin{equation}\label{ACFerror}
\Biggl[\sum_{k=0}^{N} \{ r_y(k) - r_x(k)\}^2/N\Biggr]^{1/2}
\end{equation}
for some $ N $, sufficiently large or even infinite, whe\-re~$ r_y(k)
$ and $ r_x(k) $ are the theoretical ACFs (if available) or sample
ACFs. Clearly, we can use other distances to measure the differences
of the ACFs.
(2) For highly cyclical $\{y_t\}$, we can measure the
differences between the observed and the attractor (i.e., the
limiting state) generated by the skeleton of postulated model, after
allowing for possible phase shifts. Thus, we can use the following
quasi-sample-path measure:
\begin{equation} \label{matcherror}
\min_k \sum_{t = 1}^{T} |y_{t}-x_{t+k} |/T,
\end{equation}
where $ T $ is the sample size as before.

To check the efficacy of estimation of parameters, especially in a
simulation study, we can use an obvious measure: $ \{(\hat\theta-
\theta)^\top(\hat\theta- \theta)/p\}^{1/2}$ for any estimator $
\hat\theta$ of $ \theta= (\theta_1, \ldots, \theta_p)^\top$.
Obviously, it is a~function of the number of steps $ m $ in $\operatorname{APE}(\le
m)$ or $\operatorname{AYW}(\le m)$. Note $ m = 1$ corresponds to the commonly used
estimation method based on the least squares, or the maximum
likelihood when normality is assumed. Note that the MLE is also
based on the one-step-ahead prediction for dynamical models that are
driven by Gaussian white noise. In our plotting below, results for
$\operatorname{APE}(\le1)$ and $\operatorname{AYW}(\le1)$ are not marked separately from those
for $\operatorname{APE}(\le m)$ and $\operatorname{AYW}(\le m)$ with $m
>1$.\vspace*{2pt}

\begin{figure}

\includegraphics{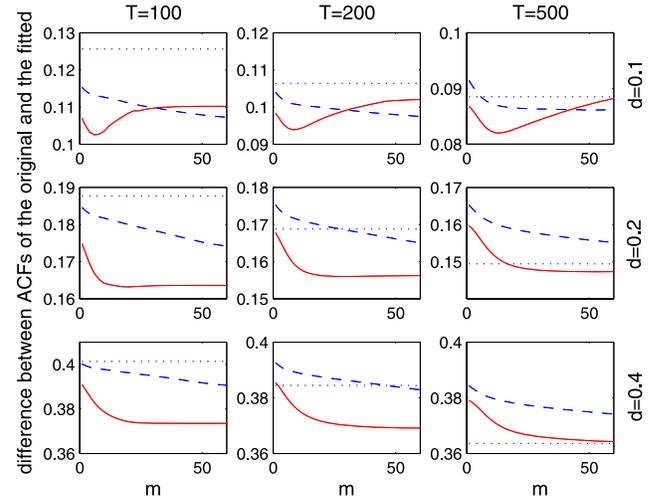}

\caption{Simulation results for Example
\protect\ref{Ex1} with different sample size $ T $, index $ d $ and the
number of steps $m$ in $\operatorname{AYW}(\le m)$ or $\operatorname{APE}(\le m)$. In each panel, the dotted
line, the solid line and the dashed line correspond to the Whittle
estimator, the APE and the AYW, respectively.} \label{figMiss}
\end{figure}

\begin{Example}[(Model misspecification)] \label{Ex1}
We postulate an $\operatorname{AR}(p)$ model to match data generated by fractionally
integrated noise $ (1-B)^d y_t = \varepsilon_t $, where $0.5
> d > -0.5$ and $ B $ is the back-shift operator and $ \{\varepsilon
_t\} $ are i.i.d. $N(0,1)$. The process is stationary,
but has long-memory when $ 0.5 > d > 0 $. The closer is~$ d $ to
$0.5$, the longer is the memory. For the use of low-order ARMA models
for short-term prediction of this type of long-memory model, see,
for example, Man (\citeyear{Man02}).
Any $\operatorname{AR}(p)$ model with finite $ p $ is a ``wrong'' model for the process.
In the following analysis, the order $p$ is
assumed unknown and determined by AIC.

The simulation results shown in Figure
\ref{figMiss} are based on 2,000 replications. We have the following
observations. (1) With a misspecified model, the $\operatorname{APE}(\le m)$ and
the $\operatorname{AYW}(\le m)$ with $m
>1$ show better matching of the ACFs than the $\operatorname{APE}(\le1)$ and
$\operatorname{AYW}(\le1)$. When $ d $ is closer to $0.5$,
the $\operatorname{AR}$ model is less likely to
fit the data well, thus necessitating a larger $m$. (2) When the
autocorrelation is not strong, which is the case with $d$ being close
to zero,
the $\operatorname{AYW}$ with large $m$ shows better matching of the ACF than the $\operatorname{APE}$;
otherwise $\operatorname{APE}$
shows better matching. It is interesting to note that although $\operatorname{APE}$ does
not target the ACF directly, it can match the ACF well in comparison
with the $\operatorname{AYW}$.
(3)~For small sample size or when $d$
is not so close to $0.5$, the $\operatorname{APE}(\le m)$
with $m >1$ show better matching than the Whittle estimator;
otherwise the Whittle estimator shows better matching.
\end{Example}

%
%

\begin{Example}[(State--space model)] \label{ssm}
$\!\!\!$Consider the $\operatorname{AR}$(4)
model with observation errors
\begin{eqnarray*}
x_t &=& \beta_1 x_{t-1} + \beta_2 x_{t-2} + \beta_3 x_{t-3} + \beta_4
x_{t-4} + \varepsilon_t,\\
  y_t &=& x_t + \eta_t.
\end{eqnarray*}
This is also a special case of a state--space model. The estimation
of the state model is of interest and has attracted considerable
attention. See, for example, Durbin and Koopman (\citeyear{DurKoo01}) and
Staudenmayer and Buonaccorsi (\citeyear{StaBuo05}).

To cover as widely as possible all admissible values on the
parameter space, we choose $ \beta_1, \beta_2, \beta_3 $ and~$\beta_4 $
uniformly distributed in the stationary region. In the
model, $ \{\varepsilon_t\} $ is a sequence of independently and
identically distributed random variables, each with a unit normal
distribution, or i.i.d. $N(0,1)$ for short; $ \{\eta_t \} $ is i.i.d. $
N(0, \sigma^2_\eta)$, such that the signal-noise ratio $
\sigma_\eta^2/\Var(y_t) = sn $ is fixed. Again, we run the
simulation 2,000 times. The results are summarized in Figures
\ref{figSS} and \ref{figSSa}. When $ p $ is known, Figure
\ref{figSS} suggests that $ \operatorname{APE}(\le m)$ and $\operatorname{AYW}(\le m)$ with $m>1$
can usually produce models that better match the dynamics of the
hidden state time series $ \{ x_t \}$ than $\operatorname{APE}(\le1)$ and
$\operatorname{AYW}(\le1)$. When $ p $ is selected by AIC, Figure~\ref{figSSa}
suggests that $\operatorname{APE}(\le m)$ and $\operatorname{AYW}(\le m)$ with $m>1$ can still
lead to better matching than $\operatorname{APE}(\le1)$ and $\operatorname{AYW}(\le1)$.

%
\begin{figure}

\includegraphics{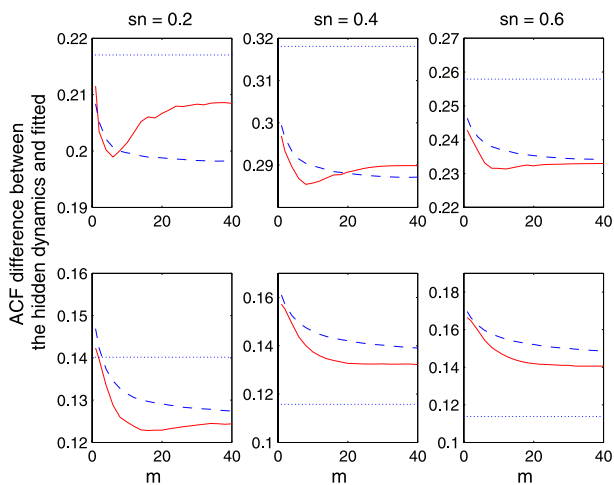}

\caption{Results for Example \protect\ref{ssm} when
the order $ p = 4 $ is known. In each panel, the dotted line, the
solid line and the dashed line correspond respectively to the Kalman
filter, the $\operatorname{APE}(\le m)$ and the $\operatorname{AYW}(\le m)$ over different
$m$.}\label{figSS}\vspace*{-3pt}
\end{figure}

\begin{figure}

\includegraphics{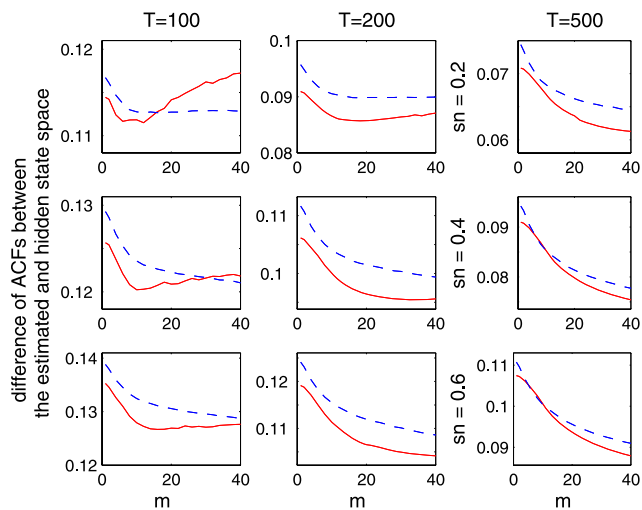}

\caption{Results for Example \protect\ref{ssm} when
the order $ p $ is selected by AIC. In each panel, the solid line is
for $\operatorname{APE}(\le m)$ and the dashed line is for $\operatorname{AYW}(\le m)$.}\label{figSSa}\vspace*{-3pt}
\end{figure}

To compare with the Kalman filter approach which utilizes the
maximum likelihood method or other methods such as the EM algorithm,
we apply the~R package ``dlm'' kindly provided by Professor Giovanni
Petris. The results are shown by dotted lines in Figure
\ref{figSS}. When the order is known, the Kalman filter shows good
performance in estimating the coefficients and in matching the ACF,
but it shows very unstable performance when the sample size is
small. Even worse, if the order is selected by the AIC, the Kalman
filter appears to be incapable of producing reasonable matching, so
much so that the results are outside the range in Figure
\ref{figSSa} in the wrong direction.
\end{Example}

\begin{Example}[(Nonlinear time series model 1:\break smooth model)] \label{Ex3}
$\!\!\!\!\!\!$Consider the simple nonlinear model
\begin{eqnarray*}
x_t &=& b_1 x_{t-1} + b_2 x_{t-1}^2 + \sigma_0 \varepsilon_t;\\
y_t &=& x_t +\sigma_1 \eta_t
\end{eqnarray*}
with parameters $ b_1 = 3.2 $ and $ b_2 = -0.2 $; both
$\varepsilon_t$ and $\eta_t$ are i.i.d. $N(0,1)$ but $\varepsilon_t$
is truncated to lie in [$-4, 4$]. We replicate our simulation 1,000
times for each set of variances $ \sigma^2_0 $ and $ \sigma^2_1 $.
The matching results are shown in Figure \ref{figMAP}.

%
\begin{figure}

\includegraphics{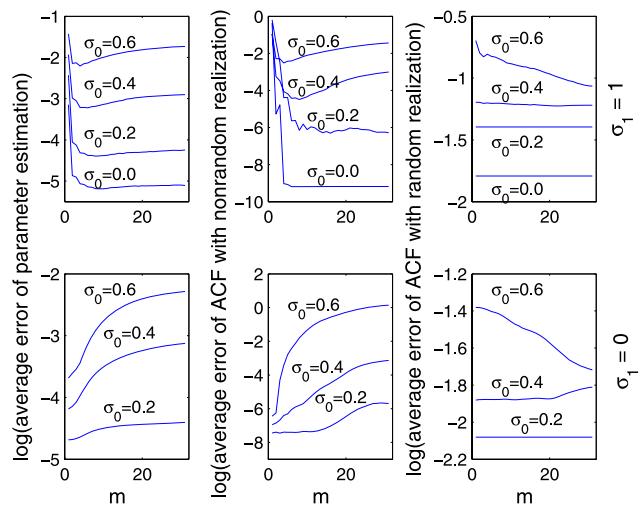}

\caption{Results for Example \protect\ref{Ex3} with $ T = 50 $ and
different $ \sigma_0 $ and $\sigma_1$. The first panel is the
estimation error of $(b_1, b_2)$ with $ \sigma_1 = 1 $; the second
panel is the difference of ACFs between the matching skeleton and
the true ACF with $ \sigma_1 = 1 $; the third panel is the
difference of ACFs between the matching skeleton and the estimated
ACFs based on random realizations with $ \sigma_1 = 1 $. Panels 4--6
are respectively the corresponding results of panels 1--3 but with $
\sigma_0 = 0 $.}\label{figMAP}\vspace*{-3pt}
\end{figure}
%
\begin{figure}

\includegraphics{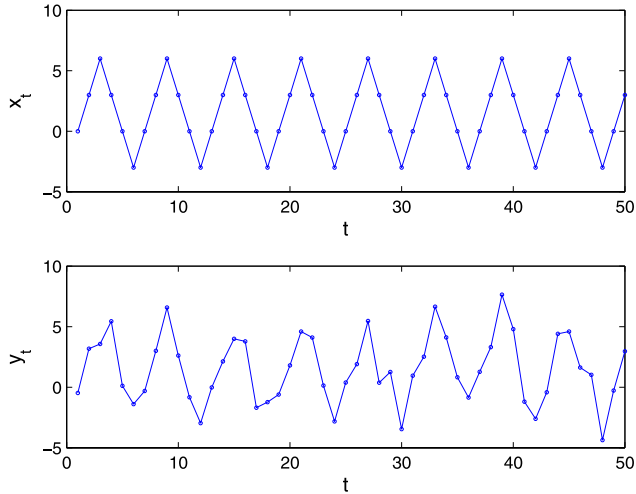}
\vspace*{-3pt}
\caption{The upper panel is a realization of
the hidden skeleton in Example \protect\ref{TARm}. The lower panel is an
observed time series subject to additive measurement error from N(0,
1).}\label{figTAR}\vspace*{-5pt}
\end{figure}

\begin{table*}[b]
\tablewidth=360pt\vspace*{-3pt}
\caption{The simulation results for Example \protect\ref{TARm}} \label{table1}
\vspace*{-3pt}\begin{tabular*}{360pt}{@{\extracolsep{\fill}}lcccc@{}}
\hline
\textbf{Model} & & \textbf{Matching} & \textbf{Cycle} & \textbf{Frequency of} \\
\textbf{setting} & \textbf{Method} & \textbf{error}& \textbf{periods} & \textbf{correct periods (\%)} \\
\hline
$T=50$, $d=2$, & $\operatorname{APE}(\le1)$ & 2.1352 (1.0334) & 5.3806 (0.6301) & 31 \\
$\mathrm{period}=6$ & $\operatorname{APE}(\le50)$ & 0.8523 (0.6591) & 5.8629 (0.5141) & 92\\
$T=50$, $d=3$, & $\operatorname{APE}(\le1)$ & 2.5301 (1.6729) & 9.4839 (0.5824) & 34 \\
$\mathrm{period}=10$ & $\operatorname{APE}(\le50)$ & 1.3987 (0.8180) & 9.9340 (0.1472) & 66\\
$T=100$, $d=2$ & $\operatorname{APE}(\le1)$ & 1.5260 (1.0643) & 5.5884 (0.6912) & 57 \\
$\mathrm{period}=6$ & $\operatorname{APE}(\le50)$ & 0.6471 (0.5301) & 5.9180 (0.3940) & 95\\
$T=100$, $d=3$ & $\operatorname{APE}(\le1)$ & 2.7196 (1.6411) & 9.4005 (0.6224) & 34\\
$\mathrm{period}=10$ & $\operatorname{APE}(\le50)$ & 1.1502 (0.5133) & 9.9705 (0.0770) & 78 \\
\hline
\end{tabular*}
\end{table*}

By coping well with noisy data due to $ \sigma_1 \eta$,\break $\operatorname{APE}(m>1)$
demonstrates substantial improvement on the parameter estimation
(in panel 1 of Figure~\ref{figMAP}), the ACF-matching of the hidden
time series\vadjust{\goodbreak} $ x_t $ (pa\-nel~2 of Figure \ref{figMAP}) and the
ACF-matching of the
observed time series (in panel 3 of Figure \ref{figMAP}). It is not
surprising that when the model is perfectly specified (i.e., $
\sigma_1 = 0)$, the $\operatorname{APE}(\le1)$ can provide better performance than
$\operatorname{APE}(\le m)$ with $m>1$ in terms of the parameter estimation and the
ACF-matching; see panels 4--5 of Figure \ref{figMAP}. However,
$\operatorname{APE}(\le m)$ with $m>1$ is still useful in matching features of the
observed time series as shown in the last panel. Our results suggest
that $ \operatorname{APE}(\le m)$ with $m>1$ leads to less improvement over $\operatorname{APE}(\le
1)$ when $\sigma_0$ (for the dynamic noise) is larger but greater
improvement when $\sigma_1$ (for the observation noise) is larger.
\end{Example}

\begin{Example}[(Nonlinear time series model 2: SETAR model)] \label{TARm}
$\!\!\!$Now, we consider a self-exciting thres\-hold autoregressive model (SETAR
model) with ske\-leton
\[
x_t = \cases{ a_0 + b_0 x_{t-1}, & if  $x_{t-d} \le c$,\cr
a_1 + b_1 x_{t-1}, & if  $x_{t-d} > c$,
}
\]\vadjust{\goodbreak}
where parameters $ a_0 = 3, b_0 = 1, a_1 = -3, b_1 = 1 $ and $ c = 0
$. A realization is shown in the first panel of Figure \ref{figTAR}.
It reveals a period of 6 when $ d = 2 $, and 10 (not shown) when $
d= 3$. Suppose that we observe
$
y_t = x_t + \eta_t,
$ where $\{ \eta_t \}$ are i.i.d. $N(0,1)$. A typical realization is
also shown in the second panel of Figure~\ref{figTAR}.

Using the APE approach to the simulated data, we denote the matching
skeleton by $ x_t $ and measure the matching error defined in 
(\ref{matcherror})
with $ T = 100 $. Based on 100 replications, we
summarize the results in Table~\ref{table1}. The matching errors have means and
standard deviations in the parentheses in column 3; the average and
standard error (in the parentheses) of the periods in all the
matching models are listed in column~4. Our results suggest that the
$\operatorname{APE}(\le m)$ with $m>1$ performs much better
than the $\operatorname{APE}(\le1)$, both in terms of matching the dynamic range and
the periodicity.
\end{Example}

\section{Application to Real Data Sets} \label{secDATA}

In this section, we study four real time series, some of which are
very well known but others less so. They are the sea levels data,
the annual sunspot numbers, Nicholson's blowflies data, and the
measles infection data in London after the massive vaccination in
the late 1960s.

\subsection{Sea Levels Data}

\begin{figure}

\includegraphics{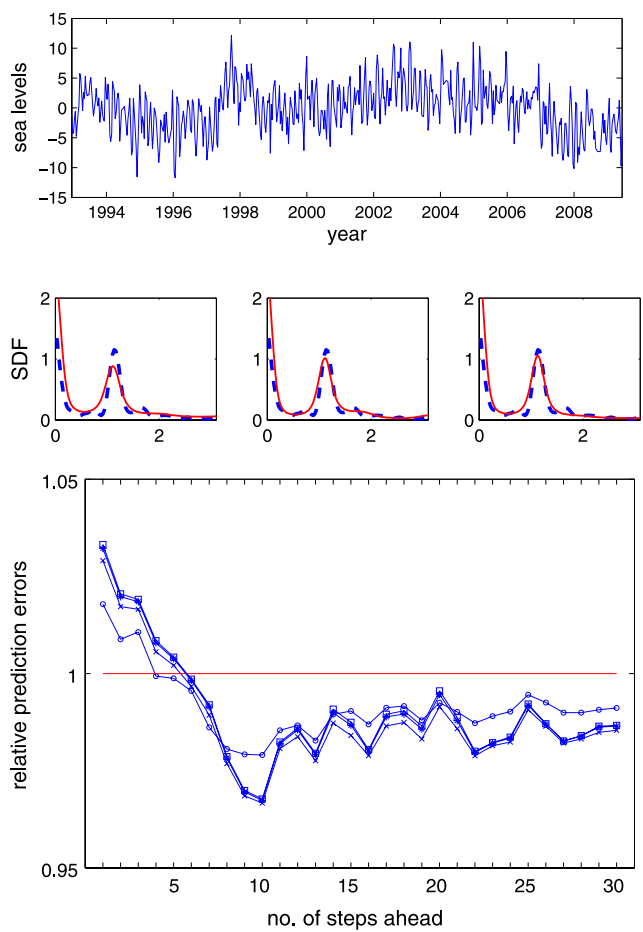}

\caption{Results for the sea level data. The data with linear
trend and seasonality removed are shown in the first panel.
Panels 2--4 are the smoothed sample SDF and those of the fitted models
by MLE,
the Whittle method and $\operatorname{APE}(\le20)$, respectively.
Panel 5 is the relative averaged multi-step-ahead
prediction errors by taking those of the one-step method as one unit.
The curves marked by ``$\circ$,'' ``$\times$,'' ``$\star$'' and ``$\diamond$'' are
for $\operatorname{APE}(\le10)$, $\operatorname{APE}(\le20)$, $\operatorname{APE}(\le30)$ and $\operatorname{APE}(\le50)$,
respectively.}
\label{figSEA}
\vspace*{-3pt}
\end{figure}

Long-term mean sea level change is of considerable interest in the
study of global climate change. Measurements of the change can
provide an important corroboration of predictions by climate models
of global warming. Starting from 1992, in each year 34 equally
spaced observations were recorded. The data with the linear trend and
seasonality removed are available at
\href{http://sealevel.colorado.edu/current/sl\_noib\_ns\_global.txt}{http://sealevel.colorado.edu/}
\href{http://sealevel.colorado.edu/current/sl\_noib\_ns\_global.txt}{current/sl\_noib\_ns\_global.txt}.
The time series is depicted in the first panel of
Figure \ref{figSEA}.
Note that the data are subject to measurement errors of 3--4 mm.

As an experiment with using a much less than ideal model to match this
data set, let us postulate an $\operatorname{AR}$ model.
By AIC, the order of the $\operatorname{AR}$ model
is selected as 6. Next, we apply the
MLE (equivalently the one-step-ahead prediction estimation method), the
Whittle method and the up-to-$m$-step-ahead
prediction estimation method to the data. The results are shown in
Figure \ref{figSEA}. The sample spectral density function (SDF) is
estimated by the method of Fan and Zhang (\citeyear{FanZha04}). The results show
clear evidence of long-memory with the singularity at the origin,
which is well captured\vadjust{\goodbreak} by all three methods. However, for the peak
away from the
origin, the Whittle estimation and $\operatorname{APE}(\le m)$ show very similar matching
capability and both show much better match than the MLE.

To investigate further, we build an $\operatorname{AR}$(6) model for every span of
observations of length $ T = 100$ and make predictions from 1 step
ahead to 30 steps ahead. For the different estimation methods, their
averaged prediction errors based on all periods are displayed in the
bottom panels of Figure \ref{figSEA}. The MLE method shows clear
superior performance for short-term prediction, while the reverse is true
from 5 steps onward.

\begin{table*}[b]
\tablewidth=400pt
\caption{The averaged difference (and its standard
error) of cycle periods in the data and
matching models and~the number of unstable matching
models [in squared brackets]} \label{table2}
\begin{tabular*}{400pt}{@{\extracolsep{\fill}}lcccc@{}}
\hline
$\bolds{m}$ \textbf{in} & \multicolumn{4}{c@{}}{\textbf{Length of time series}}
\\
\cline{2-5}
$\bolds{\operatorname{APE}(\le m)}$ & \textbf{20} & \textbf{35} & \textbf{50} & \textbf{100} \\
\hline
\phantom{0}1 & 2.5448 (3.0084) [42] & 1.7115 (1.8162) [2] & 1.3355 (1.5718) [0] &
1.5934 (1.4051) [0] \\
10 & 1.3454 (1.7082) [13] & 0.9576 (0.8499) [0] & 0.8459 (0.9584) [0] & 0.4487 (0.5427) [0] \\
20 & 1.2972 (1.7143) [10] & 0.8975 (1.1257) [0] & 0.7580 (0.6074) [0] &
0.4134 (0.9715) [0] \\
30 & & 0.8802 (1.1415) [1] & 0.8449 (0.5807) [0] & 0.3640 (0.5894) [0] \\
50 & & & 0.8548 (0.5813) [0] & 0.3538 (0.4267) [0]\\
\hline
\end{tabular*}
\end{table*}

\subsection{Annual Sunspot Numbers}

Sunspots, as an index of solar activity, are relative\-ly cooler and
darker areas on the sun's surface resul\-ting from magnetic storms.
Sunspots have a cycle of length varying from about 9 to 13 years.
\mbox{Statisticians} have fitted several models to predict sunspot numbers.
They have also noticed that the cycles are~asym\-metric
and that the time from the initial minimum~of a
cycle to its next maximum, called the rise time,~and the time from a
cycle maximum to its next \mbox{minimum}, called the fall time, are fairly
regular. Due to~their~link to other kinds of solar activity, sunspots
are helpful in predicting space weather and the state of the
ionosphere. Thus, sunspots can help predict conditions of short-wave
radio propagation as well as \mbox{satellite} communications. Historical
data of the sunspots have been recorded in different parts of the
world. The data we use are the annual sunspot numbers for the period
1700--2008 which are obtainable from
\href{http://www.ngdc.noaa.gov/stp/SOLAR/}{http://}
\href{http://www.ngdc.noaa.gov/stp/SOLAR/}{www.ngdc.noaa.gov/stp/SOLAR/}.
Yule (\citeyear{Yul27}) was the first statistician to model the sunspot number
using a model, now known as
the autoregressive model, with lag 2. Later refinements of stationary
linear models can be found in, for example, Brockwell and Davis
(\citeyear{BroDav91}) and others; higher-order $\operatorname{AR}$ models or ARMA models are used.
Akaike (\citeyear{Aka78})
suggested that the data are better modeled as nonstationary
over a long period. Tong and Lim (\citeyear{TonLim}) noticed nonlinearity in
the data dynamics and proposed the use of a self-exciting threshold
autoregressive model (or a SETAR model for short). In the following,
we postulate a two-regime SETAR model of order 3 with delay
parameter equal to 2 for the annual sunspot numbers (1700--2008).
Specifically,
\[
\hspace*{-0.09pt}x_t = \cases{ a_0 + b_0 x_{t-1} + c_0 x_{t-2} + d_0 x_{t-3}, &
if  $x_{t-2} \le\tau_0$, \cr
a_1 + b_1 x_{t-1} + c_1 x_{t-2} + d_1 x_{t-3}, & if  $x_{t-2} >
\tau_0$,
}
\]
where $ x_t = \log(\mbox{no. of sunspots} + 1) $. Note that Cheng
and Tong (\citeyear{CheTon92}) recommended\vadjust{\goodbreak} a nonparametric\break $\operatorname{AR}$(4) model. We also
tried SETAR model of order 4 with delay parameter equal to 2. The
performances of both models are very similar.

%
\begin{figure}

\includegraphics{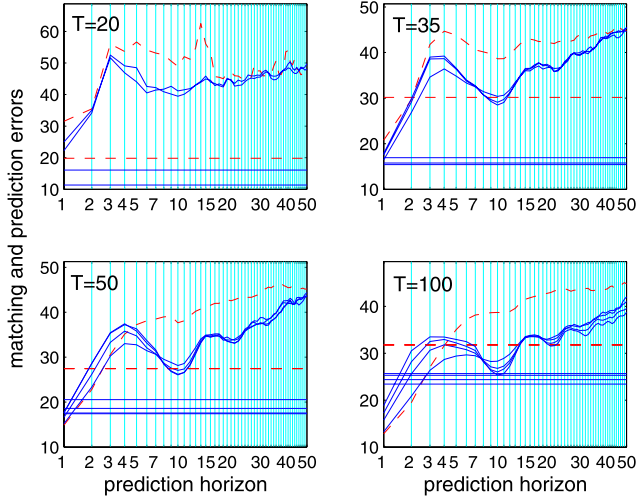}

\caption{The dashed curves are the
averaged prediction errors based on $\operatorname{APE}(\le1)$, the solid curves are
those based on $\operatorname{APE}(\le m)$ with $ m = 10, 20, 30, 50 $, respectively. The
horizontal dashed lines are the matching errors for the $\operatorname{APE}( \le1)$, the
solid lines are those for $\operatorname{APE}(\le m)$ with $ m = 10, 20, 30, 50 $,
respectively.}
\label{figsunspot}
\end{figure}

We use each fixed span of $ T $ observations to fit the postulated
model and then use it to do a post-sample prediction based on the
skeleton of the fitted model. We measure the following: (1) the
difference of cycle periods between the data and the fitted model;
(2) the frequency of stable fitted models; (3) the out-of-sample
prediction errors based on the skeletons of models fitted by the
$\operatorname{APE}(\le m)$ for different $m$; (4) the difference between the observed
time series and the time series generated by the best fitting skeleton
by reference to (\ref{matcherror}).

The results are shown in Figure \ref{figsunspot} and Table \ref{table2}. We may
draw the following conclusions.
(1)
When the observed time series is short (e.g., $T=20, 35)$, $\operatorname{APE}(\le
m)$ with $m >1$ show better matching than $\operatorname{APE}(\le1)$ in both
one-step-ahead prediction and multi-step-ahead prediction; see
panels 1 and 2 in Figure~\ref{figsunspot}. When the length of the
time series is longer (e.g., $T=50, 100)$, $\operatorname{APE}(\le1)$ can lead to
fitted models with better short-term (less than 4 steps ahead)
prediction than $\operatorname{APE}(\le m)$ with $m >1$, but for prediction beyond
4 steps ahead, the reverse appears to be the case, in line with our
understanding of the APE method. (2) When the observed time series
is short, $\operatorname{APE}(\le m)$ with $m>1$ shows its ability in avoiding
unstable models; see the numbers in the square brackets of Table~\ref{table2}.
(3) For both short time series and long time series, models fitted
by $\operatorname{APE}(\le m)$ with $m>1$ show better matching of the observed time
series in terms of their cycles; see Table \ref{table2} and the horizontal
lines in Figure~\ref{figsunspot}.

\subsection{Nicholson's Blowflies}

\begin{figure}

\includegraphics{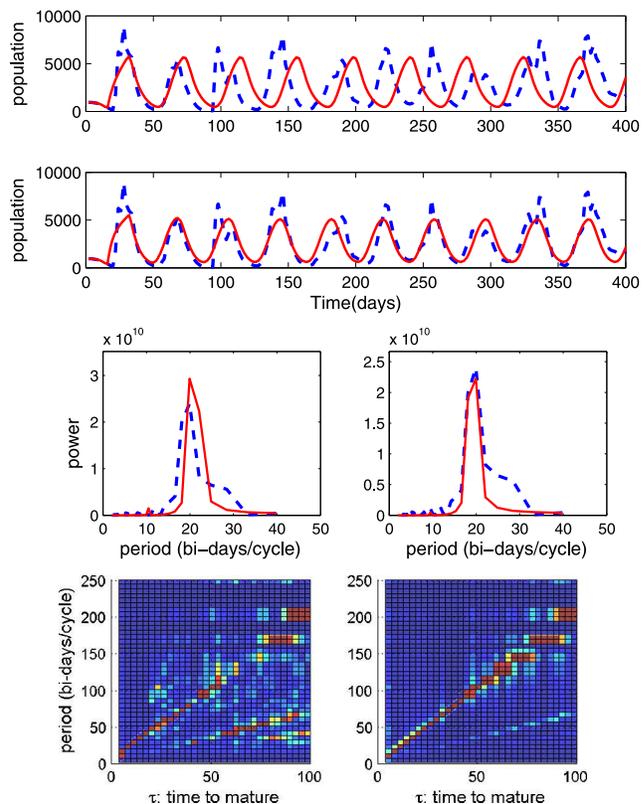}

\caption{Results for the Nicholson's blowflies data. In the
first two panels, the dashed lines are for the observed population; the
solid lines are for realizations from models fitted by $\operatorname{APE}(\le1)$
and $\operatorname{APE}(\le T)$, respectively. The dashed lines in panels 3 and 4 are the
periodograms of the observed data, and the solid
lines are those of the models fitted by $\operatorname{APE}(\le1)$ and $\operatorname{APE}(\le T)$,
respectively. In panels 5 and 6, for each $ \tau$ marked in the $x$-axis,
the vertical column is the periodogram with the values color-coded,
brighter color (blue being dull) corresponding to higher power value.
Thus the brightest point indicates the cycle-period of the dynamics
at $ \tau$.}\label{wcfig}
\end{figure}

The data consist of the total number of blowflies (Lucilia cuprina)
in a population under controlled laboratory conditions. The data
represent counts for every second day. The developmental delay
(from egg to adult) is between 14 and 15 days for the blowflies under
the conditions employed (Gurney, Blythe and
Nisbet, \citeyear{GurBlyNis80}). Nicholson
obtained 361 bi-daily recordings over a 2-year period (722 days).
However, due to biological evolution (Stokes et al., \citeyear{Stoetal88}),
the whole
series cannot be considered to represent the same system; a~major
transition appears to have occurred around day 400. Following Tong
(\citeyear{Ton90}), we consider the first part of the time series (to day 400,
thus $T=200$),
for
which the population has a 19 bi-days cycle; see Figure \ref{wcfig}.

Next, we postulate the single\vspace*{2pt} species animal population discrete
model (\ref{tpop1}) with  $ b(x_{t-\tau}) = c
x_{t-\tau}^{\alpha-1}\cdot
\exp(-x_{t-\tau}/N_0) $,
and thus
\[
x_t = c x_{t-\tau}^\alpha\exp(-x_{t-\tau}/N_0) + \nu
x_{t-1},
\]
where we take $ \tau= 8 $ (bi-days) corresponding to
the time taken for an egg to develop into an adult.
Note that we specify $ b(x_{t-\tau}) $ slightly differently from
Gurney, Blythe and
Nisbet (\citeyear{GurBlyNis80}) by adding an exponent $ \alpha- 1 $ to $x_{t-\tau
} $, which is usually necessary when a differential equation model is
discretized and approximated by a~time series model; see Glass, Xia and Grenfell (\citeyear{GlaXiaGre03}).
In the model,
there are four parameters: $ c, \alpha, N_0 $ and $ \nu$. The
(one-step-ahead prediction) MLE estimates for the parameters are
\begin{eqnarray*}
\hat c &=& 20.1192,\quad \hat N_0 = 589.5553,\\
 \hat\nu&=&
0.7598,\quad \hat\alpha= 0.8461.
\end{eqnarray*}
The APE method gives
\begin{eqnarray*}
\hat c &=& 591.5801,\quad \hat N_0 = 1307.0,\\ \hat\nu&=&
0.6469,\quad \hat\alpha= 0.2633.
\end{eqnarray*}

The skeletons based on the postulated model with parameters
estimated by above methods are shown in panels 1 and 2 in Figure
\ref{wcfig}, respectively. They show that $\operatorname{APE}(\le T)$
results in a model whose skeleton matches the observed cycles to a
much greater extent than $\operatorname{APE}(\le1)$. $\operatorname{APE}(\le1)$ gives a period
of 21 bi-days; $\operatorname{APE}(\le T)$ gives a period of 19 bi-days, which is almost
exactly the average period of the observed cycles. We have also
postulated a SETAR model. With $\operatorname{APE}(\le T)$, the SETAR model can
also capture the observed period very well, but again this is not
the case with $\operatorname{APE}(\le1)$. To investigate how the cycles change
with the time needed by the fly to grow to maturity, we vary the
time~$ \tau$ from 4 to 100 bi-days. The corresponding cycles (in
bi-days) are shown in the last two panels of Figure \ref{wcfig}.
$\operatorname{APE}(\le T)$ shows a clear linearly increasing trend in the
cycle-periods as $ \tau$ increases, while $\operatorname{APE}(\le1)$ shows
strange excursions that are difficult to interpret. The linear
relationship suggested by $\operatorname{APE}(\le T)$ may be helpful in throwing
some light on the important but not completely resolved cycle
problem of animal populations. We have also tried $\operatorname{APE}(\le m)$ with
$ m $ equal to twice or thrice the cycle-period. Their results are
similar to those of $\operatorname{APE}(\le T)$.

%


\subsection{Measles Dynamics in London}

It is well known that the continuous-time
suscepti\-ble-infected-recovered (SIR) model using a set of ordinary
differential equations can describe qualitatively the behavior of
epidemics quite well. However, it is difficult to use it for real
data modeling when the observations are made in discrete time. To
bridge the gap between the theoretical model and real data fitting,
several discrete-time or chain models have been introduced. The
Nicholson--Bailey host-parasite model (Nicholson and Bailey, \citeyear{NicBai35}) is
an early example. Bailey (\citeyear{Bai57}), Bartlett (\citeyear{BarN2}) and Finkenst\"{a}dt
and Grenfell (\citeyear{FinGre00}) proposed different types of discrete-time
epidemic models. A general discrete-time or chain model can be
written as follows:
\begin{equation}\label{tsir}
\cases{ S_{t+1} = S_t + B_t - I_{t+1}, \cr
I_{t+1} = S_t P(I_t),
}
\end{equation}
where $ I_t $, $S_t $ and $ B_t $ are respectively the number of the
infectious, the number of the susceptible and the number of births,
all at the $t$th time unit. There are many possible functional forms
for the (probability) $ P(I_t ) $. Examples are $1-(1-r_0/N)^{I_t} $
(Bartlett, \citeyear{BarN2}),
$1-\exp(-r_0I_t/N)$ (Bartlett, \citeyear{Bar56}), $r_0I_t/N $ (Baily, \citeyear{Bai57}) and
$R_0I_t^\alpha/N $ (Liu, Hethcote and Levin, \citeyear{LiuHetLev87}; Finkenst\"{a}dt and
Grenfell, \citeyear{FinGre00}), where $ N $ is the effective population of hosts,
and $ r_0 $ is the basic reproductive rate.


Next, we postulate the following (deterministic)~di\-screte-time SIR
model for
the transmission of \mbox{measles}:
\begin{eqnarray*}
I_{t+1} &=& \exp(\delta_{t,k} \beta_k) S_t I_t,\\
  S_{t+1} &=& S_{t}
+ b_t - I_{t+1} = S_0 + \sum_{\tau=0}^t B_\tau- \sum_{\tau
=1}^{t+1} I_\tau,
\end{eqnarray*}
where $ \exp(\delta_{t,k} \beta_{k}) $ is employed to indicate the
seasonality force, with $ \delta_{t,k} = 1 $ if time $ t $ is at
the $k$th season, 0 otherwise. For measles, the time unit for~$ t $
is bi-weekly, based on the infection procedure of measle; see
Finkenst\"{a}dt and Grenfell (\citeyear{FinGre00}).
Now, $k = 1, \ldots, 26$ bi-weeks
corresponds to about 54 weeks in a year. Finkenst\"{a}dt and
Grenfell (\citeyear{FinGre00}) considered the same model but with the first
equation being $ I_{t+1} = \exp(\delta_{t,k} \beta_k) S_t I_t^\alpha
$. Here, we take $ \alpha= 1 $ for two reasons. (1) If $ \alpha< 1
$, Finkenst\"{a}dt and Grenfell (\citeyear{FinGre00}) were unable to use the model to
explain the dynamics of measles in the massive vaccination era. (2)
Experience with statistical modeling of ecological populations suggests
that $ \alpha$ can be
taken as 1 with improved interpretation; see Bj{\o}nstad, Finken\-st{\"a}dt and
Grenfell (\citeyear{BjnFinGre02}).  In practice, $I_t$ may not be observed directly; what
can be observed is a random variable, say $y_t$, that has mean~$I_t$.
For this observable $y_t$, we postulate a model $x_t$ that follows a~%
Poisson distribution
with mean $I_t$.

\begin{figure}[t]

\includegraphics{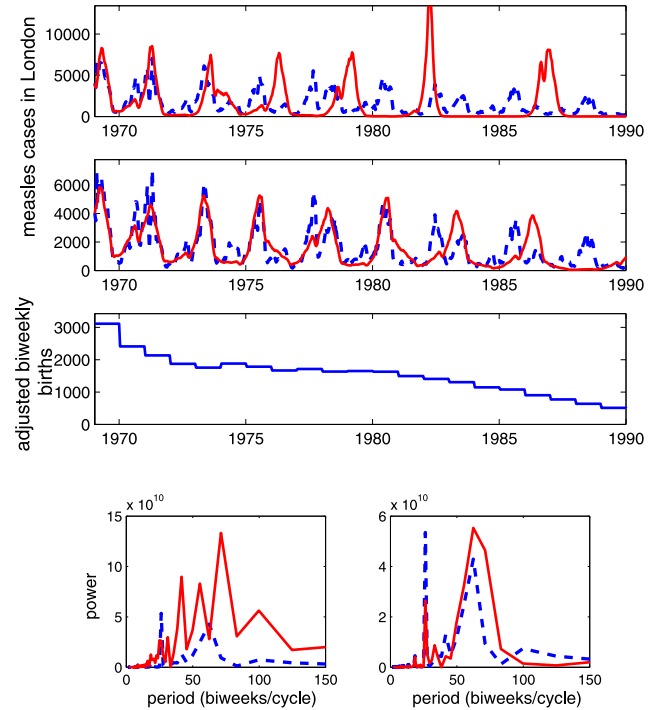}

\caption{Results for modeling the measles
incidents in London. The dashed lines in panels 1 and 2 are the
recovered incidents of measles; the solid lines are the realizations
of the model based on $\operatorname{APE}(\le1)$ and $\operatorname{APE}(\le T)$, respectively.
Panel 3 is
the adjusted birth rate by removing the vaccinated; in the bottom
panels, the dashed lines are the periodograms of the data and the red lines
are those of the matching skeleton by $\operatorname{APE}(\le1)$ and $\operatorname{APE}(\le T)$,
respectively. }\label{measlesfig}
\end{figure}

There are some problems with the data. There is nonnegligible
observation error in the data due to the under-reporting rate, which
can be as high as 50\%; see Finkenst\"{a}dt and Grenfell (\citeyear{FinGre00}), where
a~method was proposed to recover the data. Following their method,
the data were adjusted for the under-reporting rate. The adjusted
data are shown in dashed lines in panels 1 and 2 of Figure
\ref{measlesfig}. It is known that the role of vaccination is
equivalent to the reduction of the birth rate (Earn et al., \citeyear{Earetal00}).
Thus, we adjust the number of births by multiplying it by the
un-vaccination rate,
that is, $1-$(vaccination rate). We show the adjusted births in the
third panel of Figure \ref{measlesfig}. Another problem with the
data is that the susceptible $ S_t $ is unknown, which can also be
reconstructed by the method of Finkenst\"{a}dt and Grenfell (\citeyear{FinGre00}).

%

%
%
\begin{table*}
\tablewidth=400pt
\caption{Parameters in the measles transmission model} \label{tableMeasles}
\begin{tabular*}{400pt}{@{\extracolsep{\fill}}lcccccccc@{}}
\hline
\textbf{Method} & $\bolds{\beta_1}$ & $\bolds{\beta_2}$ & $\bolds{\beta_3}$ & $\bolds{\beta_4}$ & $\bolds{\beta_5}$
& $\bolds{\beta_6}$ & $\bolds{\beta_7}$ & $\bolds{\beta_8}$ \\
\hline
$\operatorname{APE}(\le1)$ &
$-$11.92&
$-$12.00&
$-$11.88&
$-$11.99&
$-$11.89&
$-$11.81&
$-$11.89&
$-$11.97 \\
$\operatorname{APE}( \le T)$ & $-$11.95&
$-$12.00&
$-$11.93&
$-$11.99&
$-$11.93&
$-$11.89&
$-$11.93&
$-$11.98
\\ [5pt]
& $\bolds{\beta_9}$ & $\bolds{\beta_{10}}$ & $\bolds{\beta_{11}}$ & $\bolds{\beta_{12}}$ & $\bolds{\beta
_{13}}$ & $\bolds{\beta_{14}}$ & $\bolds{\beta_{15}}$ & $\bolds{\beta_{16}}$ \\
$\operatorname{APE}( \le1)$ &
$-$11.92&
$-$11.99&
$-$12.05&
$-$12.01&
$-$11.93&
$-$11.96&
$-$11.98&
$-$12.04\\
$\operatorname{APE}( \le T)$ & $-$11.95&
$-$11.99&
$-$12.03&
$-$12.00&
$-$11.96&
$-$11.98&
$-$11.99&
$-$12.02
\\ [5pt]
& $\bolds{\beta_{17}}$ & $\bolds{\beta_{18}}$ & $\bolds{\beta_{19}}$ & $\bolds{\beta_{20}}$ &
$\bolds{\beta_{21}}$ & $\bolds{\beta_{22}}$ & $\bolds{\beta_{23}}$ & $\bolds{\beta_{24}}$ \\
$\operatorname{APE}( \le1)$ &
$-$11.95&
$-$12.15&
$-$12.28&
$-$12.40&
$-$12.21&
$-$11.99&
$-$11.79&
$-$11.87\\
$\operatorname{APE}( \le T)$ & $-$11.97&
$-$12.08&
$-$12.16&
$-$12.23&
$-$12.12&
$-$11.99&
$-$11.87&
$-$11.92
\\ [5pt]
& $\bolds{\beta_{25}}$ & $\bolds{\beta_{26}}$ & $\bolds{S_0}$\\
$\operatorname{APE}( \le1)$ &
$-$11.99&
$-$11.98&
17,8280 \\
$\operatorname{APE}( \le T)$ &
$-$11.99&
$-$11.98&
16,8190 \\
\hline
\end{tabular*}
\end{table*}



The estimates of the model by $\operatorname{APE}(\le1)$ are listed in Table \ref{tableMeasles}. To
ease the calculation of $\operatorname{APE}(\le T)$, we simplify the model by taking $
\beta_k = \bar\beta+ \lambda(\beta_{k,1} - \bar\beta) $, where $
\beta_{1,1}, \ldots, \beta_{26,1} $ are the estimates of $\operatorname{APE}(\le1)$
and $\bar\beta$ is their average. Consequently, only $\lambda$ and $
S_0 $ need to be estimated in implementing $\operatorname{APE}(\le T)$.
The skeletons based on models
fitted by $\operatorname{APE}(\le1)$ and $\operatorname{APE}(\le T)$ are shown in solid red lines
in panel~1 and panel 2 of Figure \ref{measlesfig}, respectively.
$\operatorname{APE}(\le T)$ shows a much better match than $\operatorname{APE}(\le1)$ in terms
of outbreak scale and cycle period. The periodogram is also much
better matched by $\operatorname{APE}(\le T)$ than by $\operatorname{APE}(\le1)$; see the last
two panels of Figure \ref{measlesfig}. We have also tried $\operatorname{APE}(\le
m)$ with $ m $ being twice or thrice the cycle period (i.e., 26
bi-weeks). The results are similar to $\operatorname{APE}(\le T)$.

\begin{figure}

\includegraphics{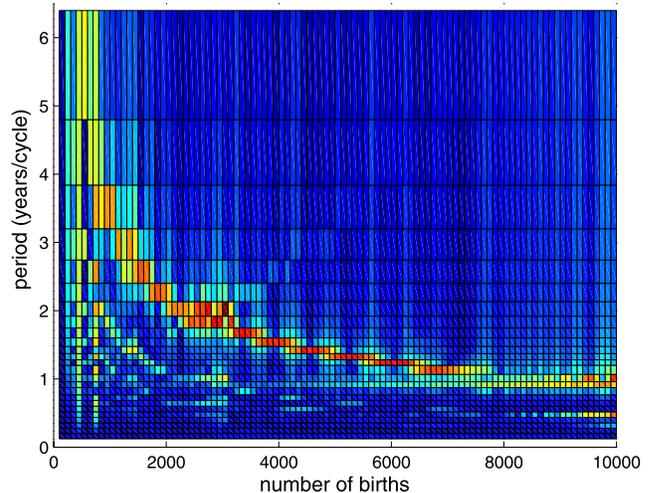}

\caption{Measles
transmission. Each vertical column is the periodogram with the
values color-coded, brighter color corresponding to higher
power. (Dark blue is considered a dull color.)}\label{measlesfig2}
\end{figure}

An important feature in the measles transmission is that there were
some big annual outbreaks in the 1950s when the birth rate was very
high after the second world war, and some big bi-annual outbreaks in
the middle of the 1960s when the birth rate was
relatively low. The
dynamics before the massive vaccination in the late 1960s was
modeled very well by a~time series model in Finkenst\"{a}dt and
Grenfell (\citeyear{FinGre00}). The theory that relates population cycle length to
birth rate has been well accepted in epidemiology and ecology.
In epidemiology, the relationship will either prolong or shorten the cumulation
procedure of susceptibles for a big outbreak.
Observations from the other sources
have lent support to this theory. For example, the measles in New
York have a three-year or four-year cycle when the birth rate is very low.
As another supporting piece of evidence, in the vaccination era, the
cycles lasted longer, to four or five years because vaccination is
equivalent to the reduction of birth rate in the transmission of
disease. However, the dynamics after the massive vaccination is
difficult to model due to the quickly changing birth rate. The
method of Finkenst\"{a}dt and Grenfell (\citeyear{FinGre00}) has failed to capture this
change of cycles in the vaccination era. It is therefore
worth noting that our
modified model, with the aid of $\operatorname{APE}(\le m)$ with $m>1$, shows satisfactory
matching. To investigate further how the cycles change with the
birth rate, for each fixed number of births we run the estimated
model and depict its periodogram and highlight the peaks by
color-coding (brighter
color for higher power). The peaks with the brightest points
correspond to the cycles of the postulated model. Figure~\ref{measlesfig2} shows clearly that when the birth rate is high
(from about 5,000 upward) the cycle is annual, but when the birth
rate is medium at about 3,000 to 4,000, the cycles become two-year
cycles. As the birth rate gets lower, the model shows that cycles
become three-year cycles or even five-year cycles.
It seems that by fitting a substantive model with the catch-all
approach, we have obtained perhaps the first \textit{discrete-time} model
that is capable of revealing the complete function linking birth-rates
to the cyclicity of measles epidemics, thereby
lending support to the general theory developed by Earn et al. (\citeyear{Earetal00}),
which was based on differential SIR equations
in continuous time.

%
%



\section{Concluding Remarks and~Further~Problems}

In this paper, we adhere to Box's dictum and abandon, right from
the very beginning, the assumption of either the postulated
parametric model being true or the observations being
error-free. Instead, we focus on ways to improve the feature
matching of
a~postulated parametric model to the observable time series.
We have introduced the notion of an optimal
parameter in the absence of a true model and defined a new form of
consistency. In particular, we have synthesized earlier
attempts into a systematic approach of estimation of the optimal
parameter, by reference to up-to-$m$-step-ahead predictions of the
postulated model.
We have also developed some general results
with proofs.

Conventional methods of estimation are typically based on just the
one-step-ahead prediction. Our~ana\-lysis, simulation study and real
applications have~con\-vinced us that they are often found wanting in
many situations, for example, the absence of a true model, short data sets,
observation errors, highly cyclical data and others.
Our stated primary objective is feature matching. Prediction is
secondary here. However, we have evidence to suggest that
a~model with good feature matching can stand a~better chance of
enjoying good medium- to
long-term prediction. Of course, if the aim is prediction with a~\textit{specified} horizon,
say $m_0$, then we simply set $w_{m_0} = 1$ and the
rest zero. In that case,
our catch-all approach really offers nothing new.

Let us now take another look at the difference between $\operatorname{APE}(\le m)$
with $m>1$ and
$\operatorname{APE}(\le1)$. Suppose we
postulate the model $ x_t = g_\theta(X_{t-1}) + \varepsilon_t $ whe\-re
$ X_{t-1} = (x_{t-1},
\ldots, x_{t-p})$ to match an observable $y$-time series. Given data
$\{y_1, y_2, \ldots, y_T\}$, $\operatorname{APE}(\le m)$ with $m>1$ and with a
constant $w_j>0$,\break all $j$,
estimates~$\theta$ by
minimizing the objective function\looseness=1
\begin{eqnarray*}
L_m(\theta) &=& \sum_{t=p+1}^T \sum_{k=1}^{\min(m, T-t)} \bigl\{
y_{t-1+k} - g_\theta^{[k]}(Y_{t-1})\bigr\}^2\\
&=& L_1(\theta)+ L_1^+(\theta),
\end{eqnarray*}\looseness=0
where $ Y_{t-1} = (y_{t-1},\ldots, y_{t-p}) $ and
\begin{eqnarray*}
L_1(\theta) &=& \sum_{t=p+1}^T \bigl\{ y_{t} - g_\theta^{[k]}(Y_{t-1})\bigr\}^2,\\
L_1^+(\theta) &=& \sum_{t=p+1}^T \sum_{k=2}^{\min(m, T-t)} \bigl\{
y_{t-1+k} - g_\theta^{[k]}(Y_{t-1})\bigr\}^2.
\end{eqnarray*}
Note that $ L_1(\theta) $ is the commonly used objective function
for $\operatorname{APE}(\le1)$, while $ L_1^+(\theta) $ is the extra information
provided by the dynamics. In terms of samples, $ L_1(\theta) $ is
based on sample $ \{ y_t, Y_{t-1}\dvtx t = p+1,\ldots, T
\} $. The extra term $ L_1^+(\theta) $ is associated with the
extra pseudo designed samples $ \{ y_{t-1+k}, Y_{t-1}\dvtx t = p+1, \ldots, T,
k = 1, \ldots, m \} $. If the data are
actually generated by the postulated model (a rare event), then
under some general conditions such as $ \varepsilon_t $ are i.i.d.
normal, $ L_1(\theta) $ will include all the information about $
\theta$. In that case, estimation based on $ L_1(\theta) $ alone is
the most efficient and the extra term $ L_1^+(\theta) $ can
provide no additional information. However, if the data are not
exactly generated by the postulated model (a common event), the extra
information provided by $ L_1^+(\theta) $ can indeed be very
helpful and should be exploited.


Despite evidence, both theoretical and practical, of the utility of
the catch-all approach, much more remains to be done. Our paper
should be seen as the first word on feature matching. Although we
have provided some concrete approaches, such as the
catch-all-conditional-mean approach, the catch-all-ACF approach,
which can easily be generalized to catch-all-$m$th-order moments and
others, there are outstanding issues. For example, we can, at
present, offer no theoretical guidance on the specification of the
weights, $\{w_m\}$. We have only offered some practical suggestions
based on our experience. It would be interesting to investigate
further possible connections with a prior in Bayesian statistics.

We have been quite fortunate
with our real examples using the APE method,
thanks to our long-standing collaboration with ecologists and epidemiologists.
However, we
are conscious of the need for the accumulation of
further experience. We are convinced that, especially in the area of
substantive modeling,
guidance by relevant subject scientists is paramount.
Relevant references include He, Ionides and King (\citeyear{HeIonKin10}), King et al. (\citeyear{Kinetal08}),
Laneri et al. (\citeyear{Lan}) and others.

Last but not least, future research should include
at least the following:
other weaker forms of
(\ref{same}), choice of a suitable weaker form in a specific application,
other criteria for model comparison, non-additive
and/or heteroscedastic measurement errors, the relaxation
of stationarity, the effect of prefiltering of data, multiple time
series, model selection among
a set of wrong models (each fitted by the catch-all method; perhaps the
idea of model calibration
in econometrics might be useful here), possible extension to other
types of
dependent data, for example, spatial data.

\appendix
\section*{Appendix: Outlines of Theoretical Justification}\label{appendix}

We need the following assumptions. However, the\-se assumptions can be
relaxed with more complicated theoretical derivation.

\begin{enumerate}[(C3)]
\item[(C1)] Time series $ \{y_t\} $ is a strictly stationary and
strongly mixing sequence with exponentially decreasing
mixing-coefficients.
\item[(C2)] The moments\vspace*{1pt} $ \EE\Vert y_{t}\Vert ^{2\delta}$, $\EE\Vert g_\vartheta
^{[k]}(y_{t}, \ldots, y_{t-p})\}\Vert ^{2\delta}$, $\EE\Vert \partial g_\vartheta
^{[k]}(y_{t}, \ldots, y_{t-p}) /\partial\theta\Vert ^\delta$ and\vspace*{2pt} $ \EE
\Vert {\partial^2
g_{\vartheta}^{[k]}(y_{t})}/\break({\partial\theta\,\partial\theta^\top
})\Vert ^\delta$ exist for some $ \delta> 2 $.

\item[(C3)] The functions $ {\partial g_{\theta
}^{[k]}(y_t)}/{\partial\theta} $ and ${\partial^2
g_{\theta}^{[k]}(y_{t})}/\break(\partial\theta\,\partial\theta^\top) $
are continuous in $ \theta\in\Theta$ and
\begin{eqnarray*}
&&\Omega\stackrel{\mathrm{def}}{=} \EE\sum_{k=1}^m w_k\biggl\{\frac{\partial
g_{\vartheta}^{[k]}(y_t)}{\partial\theta}\frac{\partial
g_{\vartheta}^{[k]}(y_t)}{\partial\theta^\top}\\
&&\hphantom{\Omega\stackrel{\mathrm{def}}{=} \EE\sum_{k=1}^m w_k\biggl\{}{}- \bigl[y_{t+k}-
g_\vartheta^{[k]}(y_t)\bigr] \frac{\partial^2
g_{\vartheta}^{[k]}(y_{t})}{\partial\theta\,\partial\theta^\top
}\biggr\}
\end{eqnarray*}
is nonsingular.

\item[(C4)] The function $ \sum_{k=1}^m w_k E [ y_{t+k} - g_\theta
^{[k]} (Y_t)]^2 $ has  a~unique minimum point for $ \theta$ in the
parameter space~$ \Theta$.
\end{enumerate}

\renewcommand{\theTheorem}{\Alph{Theorem}}
\setcounter{equation}{0}
\begin{Theorem}\label{thma}
Suppose that $\{x_t(\theta)\}$ and $\{y_t\}$
ha\-ve
the same marginal distribution and each has se\-cond-order moments.
Then
\begin{eqnarray*}
D_{_C}(y_t, x_t(\theta)) &\le& C_1 \tilde Q(y_t, x_t(\theta)),\\
D_{_F}(y_t, x_t(\theta)) &\le& C_2 \tilde Q(y_t, x_t(\theta))
\end{eqnarray*}
for some positive constants $ C_1 $ and $C_2$. Moreover, if
$\{x_t(\theta)\}$ and $\{y_t\}$ are linear $\operatorname{AR}$ models, then there are
some positive
constants
$ C_3 $ and $ C_4 $ such that
\begin{eqnarray*}
\tilde Q(y_t, x_t(\theta)) &\le& C_3 D_{_C}(y_t, x_t(\theta)),\\
\tilde Q(y_t, x_t(\theta)) &\le& C_4 D_{_F}(y_t, x_t(\theta)).
\end{eqnarray*}
\end{Theorem}

\begin{pf} By the condition on the
marginal distributions, we have
\begin{equation}\label{ee}
\E(y_{t+m}) = \E(x_{t+m}).
\end{equation}
Since $ \E[ y_{t}\{y_{t+m} - \E(y_{t+m}|y_t)\}] = 0$, we have
\begin{eqnarray*}
\E(y_{t} y_{t+m}) &=& \E\bigl\{y_{t} g^{[m]}_\theta(y_t)\bigr\} +
\E\bigl[ y_{t}\bigl\{y_{t+m} - g^{[m]}(y_t)\bigr\}\bigr]\\
&=& \E\bigl\{y_{t} g_\theta^{[m]}(y_t)\bigr\}\\
&&{} + \E\bigl[ y_{t}\bigl\{\E
(y_{t+m}|y_t) - g_\theta^{[m]}(y_t)\bigr\}\bigr].
\end{eqnarray*}
By the assumption on the marginal distribution, we have
\begin{eqnarray*}
\E\bigl\{y_{t} g^{[m]}_\theta(y_t)\bigr\}
 &=& \E\bigl\{x_{t}g_\theta^{[m]}(x_t)\bigr\}\\
 &=&\E\{x_{t}\E(x_{t+m}|x_t)\} = \E(x_{t}x_{t+m}).
\end{eqnarray*}
Thus
\begin{eqnarray}\label{dd}
\hspace*{30pt}\E(y_{t} y_{t+m}) &=& \E(x_{t} x_{t+m})\nonumber\\ [-8pt]\\ [-8pt]
&&{} + \E\bigl[ y_{t}\bigl\{\E
(y_{t+m}|y_t)- g_\theta^{[m]}(y_t)\bigr\}\bigr].\nonumber
\end{eqnarray}
It follows from (\ref{ee}) and (\ref{dd}) that
\[
\gamma_{y}(m) = \gamma_{x}(m) + \Delta_m,
\]
where
$
\Delta_m = \E[ y_{t}\{\E(y_{t+m}|y_t)- g_\theta^{[m]}(y_t)\}
]
$. By the H\"{o}lder inequality, we have
\[
|\Delta_m| \le
\{\E y^2_{t}\}^{1/2} \bigl\{\E\bigl\{\E(y_{t+m}|y_t) - g_\theta^{[m]}(y_t)\bigr\}
^2\bigr\}^{1/2}.
\]
Therefore,
\begin{eqnarray*}
&&D_{_C}(x_t(\theta), y_t)\\
&&\quad \le\sup_{\{w_k\}}\sum_{k=0}^\infty w_k \{
\E y^2_{t}\}^{1/2}\\
&&\quad\hphantom{\le\sup_{\{w_k\}}\sum_{k=0}^\infty}{}\cdot \bigl\{\E\bigl\{\E(y_{t+k}|y_t) - g_\theta^{[k]}(y_t)\bigr\}^2\bigr\}
^{1/2}\\
 &&\quad\le C_1 \tilde Q(\theta),
\end{eqnarray*}
where $ C_1 = \{\E y^2_{t}\}^{1/2}$. This is the first inequality of
Theorem \ref{thma}.

For ease of exposition, assume that $ \{y_t\} $ and $
\{x_t(\theta)\} $ are given by $\operatorname{AR}$ models with the same order, $ P $.
Otherwise we take the order as the larger of the two orders. So $ y_t =
\beta_1 y_{t-1} + \cdots+ \beta_P y_{t-P} + \varepsilon_t $ and $
x_t = \theta_1 x_{t-1} + \cdots+ \theta_P x_{t-P} + \eta_t$.

Let $ e_1 = (1, 0, \ldots, 0)^\top$, $ Y_{t-1} = (y_{t-1}, \ldots,
y_{t-P})^\top$,\break $ X_{t-1} = (x_{t-1}, \ldots, x_{t-P})^\top$, $ \mathcal{E}_t =
(\varepsilon_t, 0, \ldots, 0)^\top$ and
\begin{eqnarray*}
\Gamma_0 &=&
\pmatrix{ \beta_1 & \beta_2 & \cdots & \beta_{P-1} & \beta_P\cr
1 & 0 & \cdots & 0& 0\cr
0 & 1 & \cdots & 0& 0 \cr
\vdots& \vdots& \cdots & \vdots& \vdots\cr
0 & 0 & \cdots & 1& 0},\\
\Gamma&=&
\pmatrix{ \theta_1 & \theta_2 & \cdots & \theta_{P-1} & \theta
_{P}\cr
1 & 0 & \cdots & 0& 0\cr
0 & 1 & \cdots & 0& 0 \cr
\vdots& \vdots& \cdots & \vdots& \vdots\cr
0 & 0 & \cdots & 1 & 0
}
.
\end{eqnarray*}
Then
$
Y_{t-1+m} = e_1^\top\Gamma_0^m Y_{t-1} + e_1^\top( \mathcal{E}_{t-1+m}
+\Gamma_0 \mathcal{E}_{t-2+m} + \cdots +\Gamma_0^m \mathcal{E}_{t} ).$
It follows that
\begin{eqnarray}
&&\hspace*{10pt}[\gamma_y(m), \gamma_y(m+1), \ldots, \gamma_y(m+P-1)]\nonumber\\ [-8pt]\\ [-8pt]
&&\hspace*{10pt}\quad= \EE(y_{t-1+m}
Y_{t-1}^\top) = e_1^\top\Gamma_0^m \Sigma_0,\nonumber
\end{eqnarray}
where $ \Sigma_0 = E (Y_{t-1} Y_{t-1}^\top) = (\gamma_y(|i-j|))_{1
\le i,j \le P} $. Similarly, we have
\begin{eqnarray}
&&\hspace*{10pt}[\gamma_x(m), \gamma_x(m+1), \ldots, \gamma_x(m+P-1)]\nonumber\\ [-8pt]\\ [-8pt]
&&\hspace*{10pt}\quad = \EE(x_{t-1+m}
X_{t-1}^\top) = e_1^\top\Gamma^m \Sigma,\nonumber
\end{eqnarray}
where $ \Sigma= E (X_{t-1} X_{t-1}^\top) = (\gamma_x(|i-j|))_{1 \le
i,j \le P} $.

Assuming $ \varepsilon_t, \eta_t $ are independent sequences of
i.i.d. random variables, we have
\begin{eqnarray*}
\EE(y_{t-1+m}|Y_{t-1}) &=& e_1^\top\Gamma_0^m Y_{t-1},\\
\EE(x_{t-1+m}|X_{t-1} = Y_{t-1}) &=& e_1^\top\Gamma^m Y_{t-1}.
\end{eqnarray*}
(Note: The i.i.d. assumption can be relaxed at the expense of a much
lengthier proof.) It follows that
\begin{eqnarray*}
&& \EE\{ \EE(y_{t-1+m}|Y_{t-1}) - \EE(x_{t-1+m}|X_{t-1} = Y_{t-1})\}
^2\\
&&\quad = e_1^\top(\Gamma_0^m - \Gamma^m) \Sigma_0 (\Gamma_0^m - \Gamma
^m)^\top e_1 \\
&&\quad = e_1^\top[\Gamma_0^m \Sigma_0 - \Gamma^m \Sigma+ \Gamma^m
(\Sigma- \Sigma_0)]\\
&&\qquad{}\cdot \Sigma_0^{-1} [\Gamma_0^m \Sigma_0 - \Gamma^m
\Sigma+ \Gamma^m (\Sigma- \Sigma_0)]^\top e_1 \\
&&\quad \le\lambda_{\min}^{-1}(\Sigma_0) \|[\gamma_y(m)-\gamma
_x(m),\\
 &&\quad\hphantom{\le\lambda_{\min}^{-1}(\Sigma_0) \|[}\gamma_y(m+1)-\gamma_x(m+1), \ldots,\\
 &&\quad\hphantom{\le\lambda_{\min}^{-1}(\Sigma_0) \|[}\gamma_y(m+P-1)- \gamma
_x(m+P-1)]\\
&&\quad\hspace*{129pt}{}+ e_1^\top\Gamma^m (\Sigma- \Sigma_0) \|^2\\
&&\quad \le\lambda_{\min}^{-1}(\Sigma_0) \sum_{k=m}^{m+P-1} \{\gamma
_y(k)-\gamma_x(k)\}^2\\
&&\qquad{} + \lambda_{\min}^{-1}(\Sigma_0)\lambda
_{\max}^m(\Gamma) P \sum_{k=0}^{P-1} \{\gamma_y(k)-\gamma_x(k)\}^2,
\end{eqnarray*}
where $ \lambda_{\min}(\Sigma_0)$ and $\lambda_{\max}(\Gamma) $ are
the minimum eigenvalue of $ \Sigma_0 $ and the maximum eigenvalue
of $ \Gamma$, respectively. Note that $ \lambda_{\max}(\Gamma) < 1 $.
Therefore,
\begin{eqnarray*}
\tilde Q(\theta) &\le& P\lambda_{\min}^{-1}(\Sigma_0) \sum
_{m=0}^\infty w_m \{\gamma_y(k)-\gamma_x(k)\}^2\\
&=& C_3 D_c(x_t(\theta
), y_t),
\end{eqnarray*}
for some $ w_m \ge0 $. The proof is completed.
\end{pf}

\begin{Theorem}\label{thmb}
$\!\!\!\!$Under assumptions \textup{(C1)} and \textup{(C2)},
we have in
distribution
\[
\sqrt{n} \bigl\{ \hat\theta_{\{m\}} - \vartheta\bigr\} \to N(0, \tilde\Sigma_m),
\]
where $ \vartheta= (\tilde\Gamma_{m}^\top\tilde\Gamma_{m})^{-1}
\tilde\Gamma_{m}^\top\tilde\Upsilon_{m}$ and $ \tilde\Sigma_m $
is a positive~de\-finite matrix.
As a special case, if $ y_t = x_t + \eta_t $ with $ \Var(\varepsilon
_t)\!>\!0$ and $\Var(\eta_t)\!=\!\sigma_\eta^2
>0 $, then the above~asym\-ptotic result holds with
$ \vartheta= \theta+\sigma_\eta^2 (\Gamma_m^\top\Gamma_m + 2
\sigma_\eta^2\Gamma_p + \sigma_\eta^4 I)^{-1} (\Gamma_p +
\sigma_\eta^2 I)\theta. $
\end{Theorem}

\begin{pf} To simplify the range of
summation in the triangular array due to the lags with fixed $ m $
as $ T \to\infty$, we introduce $\cong$ to denote the fact that the
quantities on both sides of it have negligible difference.
By Theorem 3.1 of Romano and Thombs (\citeyear{RomTho96}), in an enlarged probability
space we have
\begin{eqnarray*}
\hat\Gamma_m &=& \Gamma_m + n^{-1/2} {\mathcal{U}}_m + o_p(n^{-1/2}),\\
\hat\Upsilon_m &=& \Upsilon_m + n^{-1/2} {\mathcal{V}}_m + o_p(n^{-1/2}),
\end{eqnarray*}
where $\mathcal{U}_m$ and $\mathcal{V}_m$ have the same structure as
$ \Gamma_m $ and $ \Upsilon_m $, respectively, but with $ \gamma(k) $
being replaced by $ v_k $ and $ (v_{i+1}, \ldots, v_{i+j}) $ for any
$i, j$ being jointly normal, with variance--covariance matrix given
by Romano and Thombs (\citeyear{RomTho96}). Therefore, we have
\[
\hat\theta_m = \vartheta+ n^{-1/2} \mathcal{W} + o_P(n^{-1/2}),
\]
where $\mathcal{W} = (\Gamma_m^\top\Gamma_m)^{-1} \mathcal{U}_m^\top
\Upsilon_m +(\Gamma_m^\top\Gamma_m)^{-1} \Gamma_m^\top{\mathcal{V}}_m - (\Gamma_m^\top\Gamma_m)^{-1} \{\Gamma_m^\top\mathcal{U}_m +
\mathcal{U}_m^\top\Gamma_m \} (\Gamma_m^\top\Gamma_m)^{-1} \Gamma
_m^\top$
$\Upsilon_m$ is a~linear combination of $\{v_k\}$. Thus, $\mathcal{W}$
is normally distributed with mean 0. This is the first part of Theorem~\ref{thmb}.

If $ y_t = x_t + \eta_t $, let $ \gamma_x(k) = n^{-1}
\sum_{i=1}^n x_t x_{t+k} $; it is easy to see that
\[
\hat\gamma_y(k) \cong\hat\gamma_x(k) + D_k + E_k,\quad  k = 0, 1,\ldots,
\]
where $ D_k = n^{-1} \sum_{t=1}^n (x_{t+k} + x_{t-k}) \eta_t
$ and $ E_k =\break n^{-1} \sum_{t=1}^{n} \eta_t \eta_{t+k}. $ By
the central limit theorem and Theorem 3.1 of Romano and Thombs
(\citeyear{RomTho96}), in an~en\-larged probability space there are random variables
$ \xi_k, \zeta_k $ and $ \delta_k $ such that $
\hat\gamma_x(k) = \gamma_x(k) + n^{-1/2}\xi_k + o_p(n^{-1/2}),
D_k = n^{-1/2}\zeta_k +o_p(n^{-1/2})
$ and
\[
E_k =\cases{ \sigma_\eta^2 + n^{-1/2}\delta_k + o_p(n^{-1/2}),
& if  $k =0$, \cr
n^{-1/2}\delta_k + o_p(n^{-1/2}), & if $k >0$,
}
\]
where $\xi_0, \xi_1, \ldots, \{ \zeta_k, k = 0, 1,\ldots \}, \delta_0,
\delta_1, \ldots $ are mutually independent and
$
\xi_k = \gamma_x(k) \{ \mathbf{E}\varepsilon_t^4-1\}^{1/2} W_0 + \sum
_{j=1}^\infty\{\gamma_x(j+k) + \gamma_x(j-k) \}
W_j.
$ Here $ W_0, W_1,\ldots $ are i.i.d. $N(0,1)$, $ \zeta_k \sim N(0,
2(\gamma_y(0) + \gamma_y(2k))), \Cov(\zeta_k,\allowbreak\zeta_\ell) =
2(\gamma_y(k-\ell) + \gamma_y(k+\ell)) $ and $ \delta_k \sim N(0,
\sigma_\eta^4) $ if \mbox{$ k\,{>}\,0 $} and $ \delta_0\,{\sim}\,N(0, \mathbf{E}(\eta^2\,{-}\,1)^2) $.
Define $ \Xi_k, Z_k $ and~$ \Delta_k $
similarly as $ \Gamma_k $ with $ \gamma_x(k) $ being replaced by~$
\xi_k, \zeta_k $ and $ \delta_k $, respectively. Let $ B_k$ be a $
k\times p $ matrix with the first $ p\times p $ submatrix being $
\sigma^2_\eta I_p $ and all the others 0. We have
\[
\hat\Gamma_k = \Gamma_k + B_k + n^{-1/2} \mathcal{E}_k + o_p(n^{-1/2}),
\]
where $\mathcal{E}_k = \Xi_k + Z_k + \Delta_k, $
\begin{eqnarray*}
\hat\Upsilon_k &=& \Upsilon_k + n^{-1/2} \Psi_k + o_p(n^{-1/2})\\ &=&
\Gamma_k \theta+ n^{-1/2} \Psi_k + o_p(n^{-1/2}),
\end{eqnarray*}
and $ \Psi_k = (\xi_1,\ldots, \xi_k)^\top$. It follows that
\begin{eqnarray*}
\hat\theta_m &=& [\Gamma_m ^\top\Gamma_m + 2 \sigma_\eta^2\Gamma
_p + \sigma_\eta^4 I\\
&&\hphantom{[}{} +n^{-1/2}\{ (\Gamma_m + B_m )^\top{\mathcal{E}}_m + {\mathcal{E}}_m ^\top
(\Gamma_m + B_m )\}\\
&&\hspace*{139pt}{} + o_p(n^{-1/2})]^{-1}\\
&&{}\cdot[\Gamma_m ^\top\Gamma_m + \sigma_\eta^2\Gamma_p\\
&&\hphantom{{}\cdot[}{} +n^{-1/2}\{(\Gamma_m +B_m )^\top\Psi_m + \mathcal{E}_m ^\top\Gamma_m
\theta\}\\
&&\hspace*{124pt}{} + o_p(n^{-1/2}) ] \\
&=& (\Gamma_m ^\top\Gamma_m + 2 \sigma_\eta^2\Gamma_p + \sigma
_\eta^4 I)^{-1} (\Gamma_m ^\top\Gamma_m + \sigma_\eta^2\Gamma
_p)\theta\\
&&{}+ n^{-1/2} {\mathcal{W}}_n + o(n^{-1/2}) \\
&=& \theta-
\sigma_\eta^2 (\Gamma_m ^\top\Gamma_m + 2 \sigma_\eta^2\Gamma_p
+ \sigma_\eta^4 I)^{-1} (\Gamma_p + \sigma_\eta^2 I)\theta\\
&&{}+n^{-1/2} {\mathcal{W}}_n + o(n^{-1/2}),
\end{eqnarray*}
where
$\mathcal{W}_m = (\Gamma_m ^\top\Gamma_m + 2 \sigma_\eta^2\Gamma_p
+ \sigma_\eta^4 I)^{-1} \{(\Gamma_m +\break B_m )^\top\Psi_m +\mathcal{E}_m ^\top\Gamma_m \theta\}
- (\Gamma_m ^\top\Gamma_m + 2 \sigma_\eta^2\Gamma_p + \sigma
_\eta^4 I)^{-2}\cdot\break
\{ (\Gamma_m + B_m )^\top\mathcal{E}_m
+ \mathcal{E}_m ^\top(\Gamma_m + B_m )\}
(\Gamma_m ^\top\Gamma_m + \sigma_\eta^2\Gamma_p)
$ is normally distributed. We have proved the second part.
\end{pf}

\begin{Theorem}\label{thmc}
Suppose the
system $ \{x_t = g_{\theta_0}(x_{t-1},\allowbreak\ldots, x_{t-p}) \} $ has a
finite-dimensional state--space and admits only limit cycles,
but $x_t$ is observed
as $ y_t = x_t + \eta_t $, where $ \{\eta_t \}$ are
independent with mean 0. Suppose that the function $ g_\theta(v_1,
\ldots, v_p) $ has bounded derivatives in both $ \theta$ in the
parameter space $\Theta$ and $ v_1, \ldots, v_p $ in a neighborhood
of the state--space. Suppose that the system $ z_t =
g_\theta(z_{t-1},\ldots, z_{t-p}) $ has only negative Lyapunov
exponents in a small neighborhood of $ \{x_t\} $ and in $ \theta
\in\Theta$. Let $ X_t = (x_t, x_{t-1},\ldots,\allowbreak x_{t-p}) $ and $ Y_t
= (y_t, y_{t-1}, \ldots, y_{t-p}) $.
\begin{enumerate}

\item If the observed
$Y_0 = X_0 + (\eta_{0}, \eta_{-1},\ldots, \eta_{-p}) $
is ta\-ken as the initial values of $ \{x_t\}$, then for any $n$,
\begin{eqnarray*}
&&f(y_{m+1}, \ldots, y_{m+n}|X_0 )\\
&&\quad{} - f(y_{m+1}|X_0 = Y_0)\cdots
f(y_{m+n}|X_0 = Y_0) \to0
\end{eqnarray*}
as $ m \to\infty$.

\item Suppose the equation $ \sum_{X_{t-1}} \{ g_\theta
(X_{t-1}) -x_t\}^2 = 0 $
has a unique solution in $\theta$,
where the summation is taken over all limiting states. Let $
\theta_{\{m\}} =
\operatorname{arg\,min}_{\theta} m^{-1}\!\sum_{k=1}^m\!\EE\{ $ $
y_{t-1+k}\,{-}\,g^{[k]}_\theta( Y_{t-1}) \}^2. $ If the noise takes value in
a small neighborhood of the origin, then $
\theta_{\{m\}} \to\theta_0
$
as $ m \to\infty$.
\end{enumerate}
\end{Theorem}

\begin{pf} Let $ Y_{t-1} = (y_{t-1}, \ldots, y_{t-p}), \mathcal{E}_{t-1}
=(\eta_{t-1},\allowbreak \ldots, \eta_{t-p}) $ and
$ X_{t-1} = (x_{t-1}, \ldots,  x_{t-p}) $. By the condition, we have
$ x_t = g_{\theta_0} (X_{t-1}) $. Write
\begin{eqnarray*}
&&\EE\bigl[ \bigl\{ g_\theta^{[k]} (Y_{t-1}) - x_{t-1+k}\bigr\}^2 \bigr]\\
&&\quad= \bigl\{
g^{[k]}_\theta(X_{t-1}) - g^{[k]}_{\theta_0}(X_{t-1}) \bigr\}^2 \\
&&\qquad{} - 2 \bigl\{ g^{[k]}_\theta(X_{t-1}) - g^{[k]}_{\theta_0}(X_{t-1})
\bigr\}\\
&&\hspace*{32pt}{}\cdot\EE\bigl\{ g^{[k]}_\theta( X_{t-1}+\mathcal{E}_{t-1}) - g^{[k]}_\theta(X_{t-1})\bigr\} \\
&&\qquad{} + \EE\bigl[\bigl\{ g^{[k]}_\theta( X_{t-1}+\mathcal{E}_{t-1}) - g^{[k]}_\theta
(X_{t-1})\bigr\}^2 \bigr].
\end{eqnarray*}
Note that by the definition of the Lyapunov \mbox{exponent},
\begin{eqnarray}\label{basicc}
&&\bigl| g^{[k]}_\theta( X_t+\mathcal{E}_t) - g^{[k]}_\theta(X_t)
\bigr|\nonumber\\ [-8pt]\\ [-8pt]
&&\quad \le\exp( k
\lambda) \{\EE\Vert\mathcal{E}_t\Vert\}^{1/2}.\nonumber
\end{eqnarray}
Starting from $X_0 = Y_0 $, the system at the $k$th step is~$
g^{[k]}(Y_0) $. Since the Lyapunov exponent is negative, we have
\begin{eqnarray*}
&&\bigl(g_{\theta_0}^{[m+1]}(Y_0), \ldots, g_{\theta_0}^{[m+n]}(Y_0)\bigr)\\
 &&\quad=\bigl(g_{\theta_0}^{[m+1]}(X_0), \ldots, g_{\theta_0}^{[m+n]}(X_0)\bigr)\\
 &&\qquad{} + (\delta
_{m+1}, \ldots, \delta_{m+n}),
\end{eqnarray*}
where $ \delta_{k} = g_{\theta_0}^{[k]}(Y_0) - g_{\theta
_0}^{[k]}(X_0) $,
with $ |\delta_{k} | \le\exp( k \lambda)\cdot \{\EE\Vert\mathcal{E}_0\Vert\}^{1/2} $.
Therefore,
\begin{eqnarray*}
&&(y_{m+1}, \ldots, y_{m+n})|(X_0 = Y_0)\\
&&\quad = (y_{m+1}, \ldots, y_{m+n})|X_0 +
(\delta_{m+1}, \ldots, \delta_{m+n}).
\end{eqnarray*}
Note that $ (y_{m+1}, \ldots, y_{m+n})|X_0 = (g_{\theta_0}^{[m+1]}(X_0),
\ldots,\break g_{\theta_0}^{[m+n]}(X_0)) + (\eta_{m+1}, \ldots, \eta_{m+n}) $
and that $\eta_{m+1},\ldots,\break \eta_{m+n} $ are independent.
Therefore the first part of Theorem~\ref{thmc} follows.

By (\ref{basicc}), we have
\begin{eqnarray*}
&& \bigl|\EE\bigl[ \bigl\{ g_\theta^{[k]} (Y_{t}) - x_{t+k}\bigr\}^2 \bigr] -
\bigl\{ g^{[k]}_\theta(X_t) - g^{[k]}_{\theta_0}(X_t) \bigr\}^2\bigr|\\
&&\quad \le C
\exp( k \lambda) \{\EE\Vert\mathcal{E}_t\Vert\}^{1/2}.
\end{eqnarray*}
It follows that
\begin{eqnarray*}
&& \Biggl| m^{-1} \sum_{k=1}^m \EE\bigl\{ x_{t-1+k} - g^{[k]}_\theta(Y_t)\bigr\}^2\\
&&\hphantom{\Biggl|}{} -
m^{-1} \sum_{k=1}^m \bigl\{ g^{[k]}_\theta(X_t) - g^{[k]}_{\theta_0}(X_t)\bigr\}^2 \Biggr| \\
&&\quad \le C \{\EE\Vert\mathcal{E}_t\Vert\}^{1/2} m^{-1} \sum_{k=1}^m \exp( k \lambda)\\
&&\quad \equiv\Delta(m) \to0\quad \mbox{as }  m \to\infty.
\end{eqnarray*}
That is,
\begin{eqnarray}\label{kfdg}
&&\hspace*{30pt}m^{-1} \sum_{k=1}^m \bigl\{ g^{[k]}_\theta(X_t) - g^{[k]}_{\theta_0}(X_t)
\bigr\}^2 - \Delta(m)\nonumber\\
&&\hspace*{30pt}\quad\le m^{-1} \sum_{k=1}^m \EE\bigl\{ x_{t-1+k} - g^{[k]}_\theta(Y_t) \bigr\}^2\\
&&\hspace*{30pt}\quad\le m^{-1} \sum_{k=1}^m \bigl\{ g^{[k]}_\theta(X_t) - g^{[k]}_{\theta_0}(X_t)
\bigr\}^2 +\Delta(m).\nonumber
\end{eqnarray}
By the second inequality of (\ref{kfdg}) and the continuity, we have
as $ \theta\to\theta_0 $ and $ m \to\infty$,
\begin{equation} \label{part1}
\hspace*{10pt}m^{-1} \sum_{k=1}^m \EE\bigl\{ x_{t-1+k} - g^{[k]}_\theta(Y_{t-1})\bigr\}^2
\to0.
\end{equation}
Next, we show that if $\Vert\theta- \theta_0\Vert \ge\delta>0 $, then as
$ m \to\infty$ there exists $ \delta' > 0 $ such that
\begin{equation}\label{part2}
\hspace*{25pt}m^{-1} \sum_{k=1}^m \bigl\{ g^{[k]}_\theta(X_t) - g^{[k]}_{\theta_0}(X_t)
\bigr\}^2 \ge\delta' >0.
\end{equation}
We prove (\ref{part2}) by contradiction.
Suppose the period of the limit cycle is $
\pi$. For continuous dynamics,
the assumption of a unique solution is equivalent to the statement
that as $ i \to\infty$,
\begin{eqnarray} \label{unique}
&&\sum_{k=i+1}^{i+\pi} \{ g_\theta(X_{k-1}) - x_k\}^2 \to0 \nonumber\\ [-8pt]\\ [-8pt]
&&\quad\iff\quad
\theta\to\theta_0.\nonumber
\end{eqnarray}
If (\ref{part2}) does not hold, that is, there is a $ \vartheta$ such that
\[
m^{-1} \sum_{k=1}^m \EE\bigl\{ g^{[k]}_{\vartheta} (X_t) - x_{t-1+k} \bigr\}^2
\to0,
\]
then there must be a sequence $ \{ i_j\dvtx j = 1, 2, \ldots\} $ with $ i_j
\to\infty$ as $ j \to\infty$ and
\begin{equation}\label{fad}
\hspace*{35pt}\sum_{k=i_j-p}^{i_j+\pi} \bigl\{ g^{[k]}_{ \vartheta} (X_t) - x_{t+k} \bigr\}
^2 \to0\quad \mbox{as }  j \to\infty.
\end{equation}
Let $ z_{t+k} = g^{[k]}_{ \vartheta} (X_t) $ and $ e_{t+k} = z_{t+k}-x_{t+k}$.
It follows from (\ref{fad}) that for $ k = i_j-p, \ldots, i_j+\pi$,
\[
|e_{t+k}| \to0\quad \mbox{as }  j \to\infty,
\]
and that
\begin{eqnarray*}
&&\sum_{k=i_j+1}^{i_j+\pi} \{ g_{ \vartheta}(x_{t+k-1}+e_{t+k-1}, \ldots,\\
&&\hspace*{69pt}x_{t+k-p} + e_{t+k-p})- x_{t+k} \}^2\\
&&\quad \to0.
\end{eqnarray*}
By the same argument leading to (\ref{kfdg}), we have
\begin{eqnarray*}
&&\sum_{k=i_j+1}^{i_j+\pi} \{ g_{ \vartheta}(x_{t+k-1}, \ldots,
x_{t+k-p}) - x_{t+k} \}^2\\
&&\quad \ge \sum_{k=i_j+1}^{i_j+\pi} \{ g_{
\vartheta}(x_{t+k-1}+e_{t+k-1}, \ldots,\\
&&\quad\hspace*{79pt} x_{t+k-p}+ e_{t+k-p}) - x_{t+k} \}^2\\
&&\qquad{} - C(e_{t+i_j-p}^2+\cdots+ e_{t+i_j+\pi}^2) \\
&&\quad\to 0
\end{eqnarray*}
for some $ C > 0 $. Let $ j \to\infty$; we have\break $
\sum_{k=i_j+1}^{i_j+\pi} \{ g_{ \vartheta}(x_{t+k-1}, \ldots,
x_{t+k-p}) - x_{t+k} \}^2 = 0 $, which contradicts the assumption
of
a unique solution (\ref{unique}).

By (\ref{kfdg}), (\ref{part1}) and (\ref{part2}), we have completed
the proof of Theorem \ref{thmc}.
\end{pf}

\begin{Theorem}\label{thmd}
Recall the notation in Section \ref{SecN} and
let $ \mathcal{E}_t = (\varepsilon_{t}, 0, \ldots, 0)^\top$ and $
\mathcal{N}_t = (\eta_{t}, \ldots, \eta_{t-p+1})^\top$. For
the nonlinear skeleton, we further assume that $ g_\theta(x) $ has
bounded second-order derivative with respect to $ \theta$ in
neighbor of~$ \vartheta$ for all possible values of~$ y_t $.
Suppose that the assumptions \textup{(C1)}--\textup{(C4)} hold. Then
\[
T^{-1/2} \bigl(\tilde
\theta_{\{m\}} - \vartheta_{m, \mathbf{w}} \bigr) \stackrel{D}{\to} N(0,
\Omega^{-1} \Lambda(\Omega^{-1})^\top).
\]
Specifically, for model (\ref{ssm0}) and $ y_t = x_t+ \eta_t $, if\break $
E|\varepsilon_t|^\delta< \infty$ and $ E|\eta_t|^\delta< \infty$
for some $ \delta> 4$,
then
\begin{eqnarray*}
\Lambda&=& \Cov(\Delta_{t}, \Delta_{t})\\
&&{} + \sum_{k=1}^\infty\{
\Cov(\Delta_{t}, \Delta_{t-k}) +\Cov(\Delta_{t-k}, \Delta_{t})\},\\
\Omega&=& \sum_{k=1}^m w_k\EE\biggl[ \frac{\partial g_\vartheta
^{[k]}(y_t)}{\partial\theta} \frac{\partial g_\vartheta
^{[k]}(y_t)}{\partial\theta^\top}\\
&&\hphantom{\sum_{k=1}^m w_k\EE\biggl[}{} - e_1^\top( \Phi^k - \Psi^k)X_t
\frac{\partial^2 g_\vartheta^{[k]}(X_t)}{\partial\theta\partial
\theta^\top}\\
&&\hspace*{77.7pt}{}+ e_1^\top\Phi^k \N_t \frac{\partial^2 g_\vartheta^{[k]}(\N
_t)}{\partial\theta\partial\theta^\top}\biggr]
\end{eqnarray*}
with $ \Delta_t = \sum_{k=1}^m w_k \{\sum_{j=0}^{k-1} e_1^\top
\Phi^j e_1\varepsilon_{t+k-j}\vspace*{2pt} + \eta_{t+k}\}\cdot\break {\partial
g_\vartheta^{[k]}(y_t)}/{\partial\theta}\!+\!\{e_1^\top( \Phi^k\!-\!\Psi^k)X_t\!-\!e_1^\top\Phi^k
\N_t\}\,{\cdot}\,{\partial g_\vartheta^{[k]}(y_t)}/{\partial\theta}$.
For the nonlinear model \textup{(\ref{nonlinearmodel})} and $ y_t = x_t + \eta_t $,
\begin{eqnarray*}
&&\Lambda=
\Var \Biggl[\sum_{k=1}^m w_k \frac{\partial g_\vartheta
^{[k]}(y_{t-k})}{\partial\theta} \eta_{t}\\
&&\hspace*{40.8pt}{}+ \sum_{k=1}^m w_k \bigl[g_{\theta_0}^{[k]}(X_t) - g_\vartheta
^{[k]}(y_t)\bigr] \frac{\partial g_\vartheta^{[k]}(y_t)}{\partial\theta}
\Biggr]
\end{eqnarray*}
and
\[
\Omega= \sum_{k=1}^m w_k \EE\biggl\{ \frac{\partial g_{\vartheta
}^{[k]}(y_t)}{\partial\theta}\frac{\partial g_{\vartheta
}^{[k]}(y_t)}{\partial\theta^\top} \biggr\}.
\]
\end{Theorem}

\begin{pf}
Let $ Q(\theta) = \sum_{k=1}^m w_k \EE[ y_{t+k} -
g_\theta^{[k]} (Y_t)]^2$ and
\begin{eqnarray*}
Q_n(\theta) &=& \sum_{k=1}^m w_k T^{-1} \sum_{t=1}^{T} \bigl[ y_{t+k} -
g_\theta^{[k]} (Y_t)\bigr]^2\\
 &\cong&\sum_{k=1}^m w_k \frac1{T-k} \sum
_{t=1}^{T-k} \bigl[ y_{t+k} - g_\theta^{[k]} (Y_t)\bigr]^2\\
& \stackrel
{\mathrm{def}}{=}&Q_n(\theta) .
\end{eqnarray*}
Let $ \tilde\theta_{\{m\}} = \operatorname{arg\,min}_{\theta\in\Theta}
Q_n(\theta) $. We denote this by $ \tilde\theta$ and
$\vartheta_{m,w} $ by $ \vartheta$, for simplicity. It is easy to see
that $ Q_n(\theta) \to Q(\theta) $. Following
the same argument of Wu (\citeyear{Wu1981}), we have $ \tilde\theta\to
\vartheta$ in probability.

By the definition of $ \tilde\theta$, we have $ \partial
Q_n(\tilde\theta)/\partial\theta= 0 $. By Taylor expansion, we
have
\begin{eqnarray}\label{oo1}
0 &=& \frac{\partial Q_n(\tilde\theta)}{\partial\theta}\nonumber\\ [-8pt]\\ [-8pt]
 &=& \frac
{\partial Q_n(\vartheta)}{\partial\theta} + \frac{\partial^2
Q_n(\theta^*)}{\partial
\theta\partial\theta^\top} (\tilde\theta- \vartheta),\nonumber
\end{eqnarray}
where $ \theta^* $ is a vector between $ \tilde\theta$ and $
\vartheta$, and
\begin{eqnarray}\label{firstd}
&&\hspace*{20pt}\frac{\partial Q_n(\vartheta)}{\partial\theta}\nonumber\\
&&\hspace*{20pt}\quad = - 2 T^{-1}
\sum_{k=1}^m w_k \sum_{t=1}^T \bigl[ y_{t+k} -
g_\vartheta^{[k]}(Y_t)\bigr]\nonumber\\ [-8pt]\\ [-8pt]
&&\hspace*{20pt}\quad\hphantom{= - 2 T^{-1}
\sum_{k=1}^m w_k \sum_{t=1}^T}{}\cdot\frac{\partial g_\vartheta^{[k]}(Y_t)}{\partial\theta}\nonumber\\
&&\hspace*{20pt}\quad=  2 T^{-1}\sum_{t=1}^T
\xi_{t,m},\nonumber
\end{eqnarray}
where $ \xi_{t,m} = \sum_{k=1}^m w_k [ y_{t+k} - g_\vartheta^{[k]}(Y_t)]
{\partial g_\vartheta^{[k]}(Y_t)}/{\partial\theta} $.\break
By the definition of $ \vartheta$, we have $ \partial
Q(\vartheta)/\partial\theta= 0 $, that is,
\begin{equation}\label{fkgd9}
E\xi_{t,m} = 0.
\end{equation}
Since $ y_t $ is a strongly mixing process with exponential
decreasing mixing coefficients, so is $ \xi_{t,m} $. By~(C2), we
have $ \EE\Vert\xi_{t,m}\Vert^\delta< \infty$. It follows from Theorem
2.21 of Fan and Yao [(\citeyear{FanYao2005}), page 75] that
\[
\sum_{t=1}^T \Delta_t /\sqrt{T} \stackrel{D}{\to} N\Biggl(0, \sum
_{k=0}^\infty\Gamma_\Delta(k)\Biggr).
\]
On the other hand, we have by (C3) and Proposition~2.8 of Fan and Yao
[(\citeyear{FanYao2005}), page 74]
\begin{eqnarray*}
&&\hspace*{-3pt}\frac{\partial^2 Q_n(\theta^*)}{\partial\theta\partial\theta}\\
&&\hspace*{-3pt}\quad\cong2 T^{-1} \sum_{t=1}^T \sum_{k=1}^m w_k\biggl\{\frac{\partial
g_{\theta^*}^{[k]}(Y_t)}{\partial\theta}\frac{\partial g_{\theta
^*}^{[k]}(Y_t)}{\partial\theta^\top}\\
&&\hspace*{-3pt}\quad\hphantom{\cong2 T^{-1} \sum_{t=1}^T \sum_{k=1}^m w_k\biggl\{}{}- \bigl[y_{t+k}- g_\vartheta
^{[k]}(Y_t)\bigr] \frac{\partial^2
g_{\theta^*}^{[k]}(Y_{t})}{\partial\theta\partial\theta^\top}\biggr\}\\
&&\hspace*{-3pt}\quad\to 2 \Omega.
\end{eqnarray*}

For model\vspace*{-2pt} (\ref{ssm0}), we
have $
X_{t+1} = \Phi X_t + \mathcal{E}_{t+1}
$ and
\[
X_{t+k} = \Phi^k X_t + (\mathcal{E}_{t+k} + \Phi\mathcal{E}_{t+k-1} +\cdots + \Phi
^{k-1} \mathcal{E}_{t+1}).
\]
Let $ \Psi$ be the matrix $ \Phi$ when $ \theta= \vartheta$,
respectively. Note that $ Y_t = X_t + \N_t $. It follows that
\begin{eqnarray*}
&& y_{t+k} - \Psi^k Y_t\\
&&\quad = (x_{t+k}+ \N_{t+k}) - \Phi^k
(X_{t}+ \N_{t}) + ( \Phi^k - \Psi^k)Y_k\\
&&\quad = ( \mathcal{E}_{t+k} + \Phi\mathcal{E}_{t+k-1} +\cdots + \Phi^{k-1} \mathcal{E}
_{t+1})\\
&&\qquad{} + (\N_{t+k} - \Phi^k \N_t) + ( \Phi^k - \Psi^k)Y_t
\end{eqnarray*}
and
\begin{eqnarray} \label{kfnr94}
&&y_{t+k} - e_1^\top\Psi^k Y_t\nonumber\\
&&\quad = \sum_{j=0}^{k-1} e_1^\top\Phi^j
e_1\varepsilon_{t+k-j} + \eta_{t+k}\\
&&\qquad{} + e_1^\top( \Phi^k - \Psi
^k)X_t- e_1^\top\Phi^k \N_t.\nonumber
\end{eqnarray}
It follows from (\ref{fkgd9}) and (\ref{kfnr94}) that
\[ \label{adfgd}
2 \sum_{k=1}^m w_k \EE e_1\biggl[\{( \Phi^k - \Psi^k)Y_t - \Phi
^k \N_t \} \frac{\partial g_\vartheta^{[k]}(Y_t)}{\partial\theta}\biggr]
= 0.
\]
We have
\begin{eqnarray*}
&&\hspace*{-5pt}\frac{\partial Q_n(\vartheta)}{\partial\theta}\\
&&\hspace*{-5pt}\quad \cong - 2 T^{-1}
\sum_{t=1}^T \sum_{k=1}^m w_k \Biggl[ \Biggl\{\sum_{j=0}^{k-1} e_1^\top
\Phi^j e_1\varepsilon_{t+k-j} + \eta_{t+k}\Biggr\}\\
&&\hspace*{-5pt}\hphantom{\quad \cong - 2 T^{-1}
\sum_{t=1}^T \sum_{k=1}^m w_k \Biggl[}{}\cdot \frac{\partial
g_\vartheta^{[k]}(Y_t)}{\partial\theta}\\
&&\hspace*{-5pt}\qquad{} + \Biggl\{ \{ e_1^\top( \Phi^k - \Psi^k)X_t- e_1^\top\Phi
^k \N_t\}\frac{\partial g_\vartheta^{[k]}(Y_t)}{\partial\theta} \\[-3pt]
&&\hspace*{-5pt}\qquad\hphantom{{}+\Biggl\{}{} -\EE\biggl[\{ e_1^\top( \Phi^k - \Psi^k)X_t- e_1^\top\Phi^k \N_t\}\\[-3pt]
&&\hspace*{168pt}{}\cdot\frac{\partial g_\vartheta^{[k]}(Y_t)}{\partial\theta}\biggr]\Biggr\}\Biggr]\\[-6pt]
&&\hspace*{-5pt}\quad\stackrel{\mathrm{def}}{=} - 2 T^{-1} \sum_{t=1}^T \Delta_{t}
\end{eqnarray*}
and that $ \EE\Delta_{t} = 0 $. Let $ \tilde\partial_ k = {\partial
(e_1^\top\Psi^k)}/{\partial
\theta} $. We further have
\[
\frac{\partial g_\vartheta^{[k]}(Y_t)}{\partial\theta} = \frac
{\partial e_1^\top\Psi^k}{\partial\theta} Y_t = \tilde\partial_
k (X_t + \N_t).
\]
Since $ (X_t, \varepsilon_t, \eta_t) $ is a stationary process
and a strong\-ly mixing sequence (Pham and Tran, \citeyear{PhaTra1985}) with
exponentially decreasing mixing coefficients, and $\Delta_{t}$ is a~%
function of $ \{ (X_\tau, \varepsilon_\tau, $ $ \eta_\tau)\dvtx \tau
= t, t-1, \ldots, t-m\} $, it is easy to see that $ \Delta_{t}$ is
also a strongly mixing sequence with exponentially decreasing mixing
coefficients. Note that $ \EE\Delta_{t} = 0$ and $ \EE
|\Delta_{t}|^\delta< \infty$ for some $ \delta> 2 $. By Theorem
2.21 of Fan and Yao [(\citeyear{FanYao2005}), page~75], we have
\[
\sum_{t=1}^T \Delta_t /\sqrt{T} \stackrel{D}{\to} N\Biggl(0, \sum
_{k=0}^\infty\Gamma_\Delta(k)\Biggr).
\]
On the other hand, we have in probability
\begin{eqnarray*}
&&\frac{\partial^2 Q_n(\vartheta)}{\partial\theta\partial
\theta^\top}\\[-3pt]
&&\quad= 2 T^{-1} \sum_{k=1}^m w_k \sum_{t=1}^T \frac
{\partial g_\vartheta^{[k]}(Y_t)}{\partial\theta} \frac{\partial
g_\vartheta^{[k]}(Y_t)}{\partial\theta^\top}\\
&&\qquad{} - 2 T^{-1} \sum
_{k=1}^m w_k \sum_{t=1}^T \Biggl\{ \sum_{j=0}^{k-1} e_1^\top\Phi^j
e_1\varepsilon_{t+k-j} + \eta_{t+k} \\
&&\hphantom{\qquad{} -}{} + e_1^\top( \Phi^k -
\Psi^k)X_t- e_1^\top\Phi^k \N_t\Biggr\} \frac{\partial^2
g_\vartheta^{[k]}(Y_t)}{\partial\theta\partial\theta^\top}\\
&&\quad\to 2 \sum_{k=1}^m w_k \EE\biggl\{ \frac{\partial g_\vartheta
^{[k]}(Y_t)}{\partial\theta} \frac{\partial g_\vartheta
^{[k]}(Y_t)}{\partial\theta^\top}\biggr\}\\
&&\qquad{} - 2 \sum_{k=1}^m w_k \EE
\biggl[e_1^\top( \Phi^k - \Psi^k)X_t\frac{\partial^2 g_\vartheta
^{[k]}(X_t)}{\partial\theta\partial\theta^\top}\biggr] \\
&&\qquad{} + 2\sum_{k=1}^m w_k \EE\biggl[ e_1^\top\Phi^k \N_t \frac{\partial^2
g_\vartheta^{[k]}(\N_t)}{\partial\theta\partial\theta^\top}\biggr]\\[-3pt]
&&\quad\stackrel{\mathrm{def}}= 2\Omega.
\end{eqnarray*}\vadjust{\goodbreak}
Therefore, it follows from (\ref{oo1}) that
\[
{T}^{-1/2} ( \tilde\theta-\vartheta)
\stackrel{D}{\to} N\Biggl\{0, \Omega^{-1} \sum_{k=0}^\infty
\Gamma_\Delta(k) (\Omega^{-1})^\top\Biggr\}.
\]

Next, consider model (\ref{nonlinearmodel}). Note that $
\eta_{t+k}=\break y_{t+k} - g_{\theta_0}^{[k]}(X_t) $. We have from (\ref{firstd})
that
\begin{eqnarray*}
&&\frac{\partial Q_n(\vartheta)}{\partial\theta}\\
&&\quad = - 2T^{-1} \sum_{k=1}^m w_k \sum_{t=1}^T \eta_{t+k} \frac{\partial
g_\vartheta^{[k]}(Y_t)}{\partial\theta}\\
&&\qquad - 2 T^{-1} \sum_{k=1}^m w_k
\sum_{t=1}^T \bigl[g_{\theta_0}^{[k]}(X_t) -
g_\vartheta^{[k]}(Y_t)\bigr]\\
&&\qquad\hphantom{- 2 T^{-1} \sum_{k=1}^m w_k
\sum_{t=1}^T}{}\cdot\frac{\partial g_\vartheta^{[k]}(Y_t)}{\partial\theta}\\
&&\quad\cong -2 T^{-1} \sum_{t=1}^T \Biggl\{ \Biggl[\sum_{k=1}^m w_k
\frac{\partial g_\vartheta^{[k]}(Y_{t-k})}{\partial\theta}
\Biggr]\eta_{t}\\
&&\qquad\hphantom{-2 T^{-1} \sum_{t=1}^T \Biggl\{{}}{}+ \sum_{k=1}^m w_k \bigl[g_{\theta_0}^{[k]}(X_t) - g_\vartheta
^{[k]}(Y_t)\bigr]\\
&&\hspace*{145.7pt}\qquad\cdot \frac{\partial g_\vartheta^{[k]}(Y_t)}{\partial\theta}
\Biggr\}.
\end{eqnarray*}
Let
\begin{eqnarray*}
C_m(x_{t-k}, \eta_{t-k}) &=& \Biggl[\sum_{k=1}^m w_k \frac{\partial
g_\vartheta^{[k]}(Y_{t-k})}{\partial\theta} \Biggr],\\
B_{m}(x_t, \eta_t) &=& \sum_{k=1}^m w_k \bigl[g_{\theta_0}^{[k]}(X_t) - g_\vartheta
^{[k]}(Y_t)\bigr]\\
&&\hphantom{\sum_{k=1}^m}\cdot \frac{\partial g_\vartheta^{[k]}(Y_t)}{\partial\theta}.
\end{eqnarray*}
By (\ref{fkgd9}), we have $\EE B_{m}(X_t, \eta_t)\,{=}\,0 $. Thus $
B_{m}(x_t, \eta_t) $ are independent with expectation 0. It is
easy to see that $ \xi_{m,t} =C_m(X_{t-k},
\eta_{t-k})\eta_{t} +B_{m}(X_t, \eta_t) $ is a~%
martingale difference. The Lyapunov's condition is satisfied. Thus,
we have
\begin{equation}\label{kdfjkg}
{T}^{-1/2}\sum_{t=1}^T \xi_{t, m} \stackrel{D}{\to} N\{0, \EE(\xi
_{m,t} \xi_{m,t}^\top)\}.
\end{equation}
Similarly to $ {\partial Q_n(\vartheta)}/{\partial\theta} $ above,
we have
\begin{eqnarray}\label{afitg}
&&\hspace*{20pt}\frac{\partial^2 Q_n(\theta^*)}{\partial\theta\partial
\theta}\nonumber\\
&&\hspace*{20pt}\quad\cong -2 T^{-1} \sum_{t=1}^T \Biggl[\sum_{k=1}^m w_k \frac{\partial^2
g_{\theta^*}^{[k]}(Y_{t-k})}{\partial\theta\partial\theta^\top}
\Biggr]\eta_{t}\nonumber\\ [-8pt]\\ [-8pt]
&&\hspace*{20pt}\qquad{}+ 2 T^{-1} \sum_{t=1}^T \sum_{k=1}^m w_k \frac{\partial g_{\theta
^*}^{[k]}(Y_t)}{\partial\theta}\frac{\partial g_{\theta
^*}^{[k]}(Y_t)}{\partial\theta^\top}\nonumber \\
&&\hspace*{20pt}\quad\to 2 \sum_{k=1}^m w_k \EE\biggl\{ \frac{\partial g_{\vartheta
}^{[k]}(Y_t)}{\partial\theta}\frac{\partial g_{\vartheta
}^{[k]}(Y_t)}{\partial\theta^\top} \biggr\}\\
&&\hspace*{20pt}\quad\stackrel{\mathrm{def}}= 2 \Omega.\nonumber
\end{eqnarray}
Finally, from (\ref{oo1}), (\ref{kdfjkg}) and (\ref{afitg}) we have
\[
T^{-1/2} ( \tilde\theta- \vartheta) \stackrel{D}\to N\{0, \Omega
^{-1}\EE(\xi_{m,t} \xi_{m,t}^\top) \Omega^{-1} \}.
\]
We have completed the proof.
\end{pf}

\section*{Acknowledgments} Yingcun Xia's research is supported in
part by a~grant from the Risk Management Institute, National University of
Singapore. Howell Tong gratefully acknowledges partial support from the National
University of Singapore (Saw Swee Hock Professorship) and the
University of Hong Kong (Distingui\-shed Visiting Professorship). We are
grateful to the Executive
Editor and two anonymous referees for constructive comments. We are
also grateful to the Institute of Mathematical Science, National
University of Singapore, for giving us the opportunity to present our
work at
their Workshop on Nonlinear Time Series Analysis in February, 2011.

%

\end{document}